\documentclass[preprint,3p]{elsarticle}

\usepackage{graphicx}
\usepackage{amsmath}
\usepackage{array}
\usepackage{amssymb}
\usepackage{multirow}
\usepackage{subfig}
\usepackage{color}

\newcommand{\pa}{\partial}
\newcommand{\na}{\boldsymbol{\nabla}}
\newcommand{\lam}{\lambda}

\newcommand{\ClVarr}{black}

\journal{Elsevier}

\begin{document}

\begin{frontmatter}



\title{Exact Charge-Conserving Scatter-Gather Algorithm for Particle-in-Cell Simulations on Unstructured Grids: A Geometric Perspective}


\author[rvt]{Haksu Moon\corref{cor1}}
\ead{moon.173@osu.edu}

\author[rvt]{Fernando L. Teixeira}
\ead{teixeira@ece.osu.edu}

\author[rvt2]{Yuri A. Omelchenko}
\ead{omelche@gmail.com}

\cortext[cor1]{Corresponding author}

\address[rvt]{ElectroScience Laboratory, The Ohio State University, Columbus, OH 43212, USA}

\address[rvt2]{Trinum Research Inc., San Diego, CA 92126, USA}

\begin{abstract}
We describe a charge-conserving scatter-gather algorithm for particle-in-cell simulations on unstructured grids. Charge conservation is obtained from first principles, i.e., without the need for any post-processing or correction steps. 
This algorithm recovers, at a fundamental level, the scatter-gather algorithms presented recently by Campos-Pinto et al.~\cite{Pinto14:Charge} (to first-order) and by Squire et al.~\cite{Squire12:Geometric}, but it is derived here in a streamlined fashion from a geometric viewpoint. Some ingredients reflecting this viewpoint
are (1) the use of (discrete) differential forms of various degrees to represent fields, currents, and charged particles and provide localization rules for the degrees of freedom thereof on the various grid elements (nodes, edges, facets), (2) use of Whitney forms as basic interpolants from discrete differential forms to continuum space, and (3) use of a Galerkin formula for the discrete Hodge star operators (i.e., ``mass matrices'' incorporating the metric datum of the grid) applicable to generally irregular, unstructured grids. The expressions obtained for the scatter charges and scatter currents are very concise and do not involve numerical quadrature rules. Appropriate fractional areas within each grid element are identified that represent  scatter charges and scatter currents within the element, and a simple geometric representation for the (exact) charge conservation mechanism is obtained by such identification.
The field update is based on the coupled first-order Maxwell's curl equations to avoid spurious modes with secular growth (otherwise present in formulations that discretize the second-order wave equation). 
Examples are provided to verify preservation of
discrete Gauss' law for all times.
\end{abstract}


\end{frontmatter}

% \linenumbers

% main text

%%%%%%%%%%%%%%%%%%%%%%%%%%%%%%%%%%%%%%%%%%%%%%%%%%%%%%%%%%%%%%%%%%%%%%%%%%%%%%%%%%%%%%%%%%%%%%%%%%%%%%%%%%%%%%
\section{Introduction}
%%%%%%%%%%%%%%%%%%%%%%%%%%%%%%%%%%%%%%%%%%%%%%%%%%%%%%%%%%%%%%%%%%%%%%%%%%%%%%%%%%%%%%%%%%%%%%%%%%%%%%%%%%%%%%
Particle-in-cell (PIC) algorithms have been extensively used for several decades in the simulation of problems involving space charges \cite{Hockney:Computer, Birdsall:Plasma, Fehske:Computational, Dawson83:Particle}, including plasma-related applications such as electron accelerators \cite{Bruhwiler01:Particle}, laser ignited devices \cite{Strozzi12:Fast}, and high-power microwave generation \cite{Booske08:Plasma}, to name a few.
A key challenge to PIC algorithms is how to \textcolor{\ClVarr}{achieve} exact charge conservation properties on unstructured, irregular grids. 
The traditional approach to enforce charge conservation is to apply correction terms to the electric fields in order to satisfy Gauss' law at every time step \cite{Marder87:Method, Langdon92:Enforcing, Mardahl97:Charge}.
An alternative approach is to enforce the (discrete) continuity equation directly at every grid cell~\cite{Eastwood91:Virtual, Villasenor92:Rigorous, Esirkepov01:Exact, Umeda03:New, Kong11:Particle, Sokolov13:Alternating}. However, this latter approach is predicated on the use of rectangular grids. In order to more accurately represent general curved geometries, unstructured grids are highly desirable \cite{Squire12:Geometric, Jacobs06:High}.

A charge-conserving scatter-gather algorithm for irregular grids based on first principles, that is without resorting to any correction steps, was proposed in~\cite{Candel07:Parallel, Candel09:Parallel}. This algorithm relies on the vector-wave equation and on the use of the time-integrated electric field as a dynamical variable. Compared to Maxwell's equations, the vector-wave equation admits an enlarged solution space that includes gradient-like solutions exhibiting secular growth in time, i.e., of the form $t \, \na \phi$. These solutions, even if not initially excited by (properly set) initial conditions, can nevertheless emerge at late times due to the accumulation of round-off errors and pollute the numerical solution unless specialized strategies such tree-cotree decomposition (gauging)~\cite{Manges95:cotree},
grad-div regularization~\cite{Albanese88:eddy}, or ad hoc corrections~\cite{Hwang99:stability} are utilized. In addition, the approach in~\cite{Candel07:Parallel, Candel09:Parallel} requires a numerical differentiation in time to compute the electric field $\mathbf{E}$. This causes the (temporal) order of accuracy for $\mathbf{E}$ to be one order less than the order of accuracy of the time integration scheme itself. Further, a Newmark-beta scheme is adopted in \cite{Candel07:Parallel, Candel09:Parallel} for the numerical time integration. This scheme has the advantage of producing an unconditionally stable update, but has the disadvantage of yielding a linear system with deteriorating condition numbers for large Courant factors that may occur, for example, in highly refined grids or multiscale problems~\cite{Moon14:Trade}.

Another exact charge-conserving algorithm based on first-principles was recently presented by Squire et al.~\cite{Squire12:Geometric}. 
This algorithm utilizes a variational vector-potential formulation that is multi-symplectic and has manifest gauge symmetry. The authors employ discrete Hodge star operators (``mass matrices,'' which encode the spatial metric) represented as {\it diagonal} matrices. This diagonal representation is only adequate for Delaunay triangular (primal) grids, wherein a particular type of dual grids can be constructed such that the paired primal-dual grid elements are {\it orthogonal} to each other (constituting the so-called Voronoi dual). This diagonal representation is not suited for more general types of triangular grids where a dual grid with such orthogonality property may not exist.

More recently, Campos Pinto et al.~\cite{Pinto14:Charge} put forth a comprehensive charge-conserving PIC algorithm based on a finite element Maxwell solver using curl-conforming elements of arbitrary orders, arbitrary shape factors, and piecewise polynomial trajectories of particles.  
In this paper, we derive a charge-conserving scatter-gather scheme for PIC simulations on unstructured grid that recovers, at a fundamental level, the scatter-gather algorithm by Campos Pinto et al.~\cite{Pinto14:Charge} to first-order, and by Squire et al.~\cite{Squire12:Geometric}, but is obtained
here in a more streamlined fashion from a geometric viewpoint.
Similarly to ~\cite{Pinto14:Charge} but in contrast to~\cite{Squire12:Geometric}, the present algorithm does not require a Delaunay triangular grid, being equally applicable
to any irregular triangular (or simplicial\footnote{Recall that a simplicial grid is such that all its elements are {\it simplices}, i.e., elements whose boundaries are the union of a minimal number of lower-dimensional elements. Therefore, in a 3-D simplicial grid for example, any face is a triangle and any volume cell is a tetrahedron.}) grid through the use of a sparse, but nondiagonal representation for the discrete Hodge star operators~\cite{Kim11:Parallel, He06:Geometric}, as given by expressions \eqref{hodge1} and \eqref{hodge2} below. 
The present algorithm uses the coupled first-order Maxwell's curl equations and a mixed Whitney form representation for the electromagnetic fields to avoid spurious solutions, and a leap-frog time update that only requires a symmetric-positive-definite linear system solver with no condition number deterioration across different mesh-refinement scales~\cite{Moon14:Trade}.
Some of the ingredients reflecting the geometric perspective adopted here are: (1) the use of (discrete) differential forms of various degrees to represent all dynamic variables (fields, currents, and particles) and provide clear localization rules for the degrees of freedom thereof on the appropriate grid elements (nodes, edges, facets), and (2) the use of  0, 1, and 2 (or nodal, edge, and face) Whitney  forms (interpolatory functions) to consistently represent the particle charges, currents, and fields in continuum space~\cite{Kim11:Parallel,Bossavit88:Whitney,Bossavit02:Generating}. 
In particular, the expressions obtained for the scatter charges and scatter currents are very simple and do not involve any numerical quadrature rules. Appropriate fractional areas within each grid element are identified to represent  scatter charges and scatter currents, and a geometric demonstration of the exact charge conservation is obtained from this very identification and 
 irrespective of triangular shape of the grid cells.

% end red-colored text

%%%%%%%%%%%%%%%%%%%%%%%%%%%%%%%%%%%%%%%%%%%%%%%%%%%%%%%%%%%%%%%%%%%%%%%%%%%%%%%%%%%%%%%%%%%%%%%%%%%%%%%%%%%%%%
\section{Formulation}
%%%%%%%%%%%%%%%%%%%%%%%%%%%%%%%%%%%%%%%%%%%%%%%%%%%%%%%%%%%%%%%%%%%%%%%%%%%%%%%%%%%%%%%%%%%%%%%%%%%%%%%%%%%%%%
%%%%%%%%%%%%%%%%%%%%%%%%%%%%%%%%%%%%%%%%%%%%%%%%%%%%%%%%%%%%%%%%%%%%%%%%%%%%%%%%%%%%%%%%%%%%%%%%%%%%%%%%%%%%%%
\subsection{Field update}
\label{field.update}
%%%%%%%%%%%%%%%%%%%%%%%%%%%%%%%%%%%%%%%%%%%%%%%%%%%%%%%%%%%%%%%%%%%%%%%%%%%%%%%%%%%%%%%%%%%%%%%%%%%%%%%%%%%%%%
On unstructured grids, the electric field intensity $\mathbf{E}(\mathbf{r},t)$ and the magnetic flux density $\mathbf{B}(\mathbf{r},t)$ can be expanded using Whitney basis functions as~\cite{Bossavit88:Whitney, Bossavit02:Generating,Kim11:Parallel}
\begin{flalign}
\mathbf{E}(\mathbf{r},t) = \sum^{N_e}_{i=1} e_i(t) \mathbf{W}_i^1(\mathbf{r}), \label{E}\\
\mathbf{B}(\mathbf{r},t) = \sum^{N_f}_{i=1} b_i(t) \mathbf{W}_i^2(\mathbf{r}), \label{B}
\end{flalign}
where $N_e$ and $N_f$ are the number of edges and faces in the grid, so that there is a 1:1 correspondence (localization) of the degrees of freedom  $e_i(t)$ and $b_i(t)$ to edges and faces, resp., in the grid.
In the above, $\mathbf{W}_i^1(\mathbf{r})$ and $\mathbf{W}_i^2(\mathbf{r})$ are (Whitney) edge and face basis functions \cite{Bossavit88:Whitney, Bossavit02:Generating,Kim11:Parallel}, respectively. The units of $\mathbf{W}_i^1(\mathbf{r})$ and $\mathbf{W}_i^2(\mathbf{r})$ are [m$^{-1}$] and [m$^{-2}$], respectively. The edge and face Whitney functions above are vector proxies of Whitney 1-forms and Whitney 2-forms, respectively. For details about Whitney functions, see \ref{app.a}. Note that the units of the $e_i(t)$ and $b_i(t)$ factors are Volts [V] and Webers [Wb], respectively. The above expansions are informed by the fact that electric field is a 1-form and magnetic flux density is a 2-form in the language of differential forms~\cite{Deschamps81:Electromagnetics, Warnick97:Teaching, Teixeira99:Lattice, He07:Differential}. Furthermore, if an electric current density is present in the grid, current density is defined such that
\begin{flalign}
\mathbf{J}_\star(\mathbf{r},t) = \sum^{N_e}_{i=1} i_i(t) \mathbf{W}_i^1(\mathbf{r}). \label{J}
\end{flalign}
so that the degrees of freedom $i_i(t)$ are associated to the edges of the grid, like those of $\mathbf{E}(\mathbf{r},t)$~\footnote{We employ a star subscript on $\mathbf{J}_\star$ because, strictly speaking, the electric current density $\mathbf{J}$  is a (twisted) 2-form that should be discretized in the dual grid. The above $\mathbf{J}_\star$ is the Hodge dual representation of $\mathbf{J}$, on the primal grid \cite{Teixeira99:Lattice}.}.
With the aid of Galerkin testing, Maxwell's equations can be spatially discretized as~\cite{Kim11:Parallel}
\begin{flalign}
\mathbf{C} \cdot \mathbf{e} &= -\frac{d}{dt}\mathbf{b}, \label{semi.dis.curl.1}\\
\mathbf{C}^T \cdot \left[\star_{\mu^{-1}}\right] \cdot \mathbf{b}
	&= \frac{d}{dt}\left[\star_{\epsilon}\right] \cdot \mathbf{e} + \mathbf{i}. \label{semi.dis.curl.2}
\end{flalign}
$\mathbf{C}$ is an incidence matrix with elements in the set \{-1,0,1\}, providing the (discrete) representation of curl operator distilled from the metric \cite{Clemens01:Discrete, Schuhmann01:Conservation}. The superscript $T$ indicates the transpose. The arrays of degrees of freedom are defined as
\begin{flalign}
\mathbf{e}&=[e_1(t),e_2(t),\cdots ,e_{N_e}(t)]^T, \\
\mathbf{b}&=[b_1(t),b_2(t),\cdots ,b_{N_f}(t)]^T, \\
\mathbf{i}&=[i_1(t),i_2(t),\cdots ,i_{N_e}(t)]^T.
\end{flalign}
In addition, $\left[\star_{\mu^{-1}}\right]$ and $\left[\star_{\epsilon}\right]$ in \eqref{semi.dis.curl.2} are discrete Hodge star operators given by the following integrals~\cite{Kim11:Parallel, He06:Geometric},
\begin{eqnarray}
\left[ \star_{\epsilon} \right]_{ij} & = & \int_{\Omega} \epsilon \,
\mathbf{W}_i^1(\mathbf{r}) \cdot \mathbf{W}_j^1(\mathbf{r}) \, dV, \label{hodge1} \\
\left[ \star_{\mu^{-1}} \right]_{ij} & = & \int_{\Omega} \frac{1}{\mu} \, \mathbf{W}_i^2(\mathbf{r}) \cdot \mathbf{W}_j^2(\mathbf{r}) \, dV, \label{hodge2}
\end{eqnarray}
which, for a given grid, are pre-computed using quadratures. Both $\left[\star_{\mu^{-1}}\right]$ and $\left[\star_{\epsilon}\right]$ are symmetric positive-definite matrices, which is a property is necessary to ensure stability of the time updating scheme~\cite{Teixeira99:Lattice}. It should be pointed out that the Hodge matrix associated with $\mathbf{i}$ in \eqref{semi.dis.curl.2} is the identity matrix.

Using the leap-frog scheme, the semi-discrete equations \eqref{semi.dis.curl.1} and \eqref{semi.dis.curl.2} can be fully discretized as
\begin{flalign}
\mathbf{b}^{n+\frac{1}{2}}
	&=\mathbf{b}^{n-\frac{1}{2}} - \Delta t \, \mathbf{C}\cdot\mathbf{e}^n, \label{b.update} \\
\left[\star_{\epsilon}\right] \cdot \mathbf{e}^{n+1}
	&=\left[\star_{\epsilon}\right] \cdot \mathbf{e}^n + \Delta t
	\left(
 \mathbf{C}^T \cdot \left[\star_{\mu^{-1}}\right] \cdot \mathbf{b}^{n+\frac{1}{2}} - \mathbf{i}^{n+\frac{1}{2}}
	\right). \label{e.update}
\end{flalign}
Since \eqref{hodge1} and \eqref{hodge2} are positive-definite, it can be easily be shown that the above update scheme is conditionally stable, obeying a Courant-like stability criterion~\cite{Kim11:Parallel},\cite{Schuhmann01:Conservation}. From the discrete values $\mathbf{e}^{n+1}$ and $\mathbf{b}^{n+\frac{1}{2}}$ obtained from \eqref{b.update} and \eqref{e.update},
the temporal coefficients $e_i(t)$ in \eqref{E} and $b_i(t)$  in \eqref{B} can be interpolated as~\cite{Lee06:Note}.
\begin{flalign}
e_i(t) &= \sum_n e_i^n\Pi^n(t), \label{e.time.exp}\\
b_i(t) &= \sum_n b_i^{n+\frac{1}{2}}\Lambda^{n+\frac{1}{2}}(t), \label{b.time.exp}
\end{flalign}
where $\Pi^n(t)$ is a piecewise constant (pulse) function centered on integer times and $\Lambda^{n+\frac{1}{2}}(t)$ is a piecewise linear (rooftop) function centered on half-integer times.
The choice of \eqref{e.time.exp} and \eqref{b.time.exp} is inspired by \eqref{semi.dis.curl.1}, where time derivative of $\mathbf{b}$ should be the same form of $\mathbf{e}$. Also, $i_i(t)$ in \eqref{J} can be likewise expanded in the piecewise linear function centered on half-integer time indices (see \eqref{semi.dis.curl.2}).
\begin{figure}[t]
	\centering
	\subfloat[\label{tophat}]{%
      \includegraphics[width=2.9in]{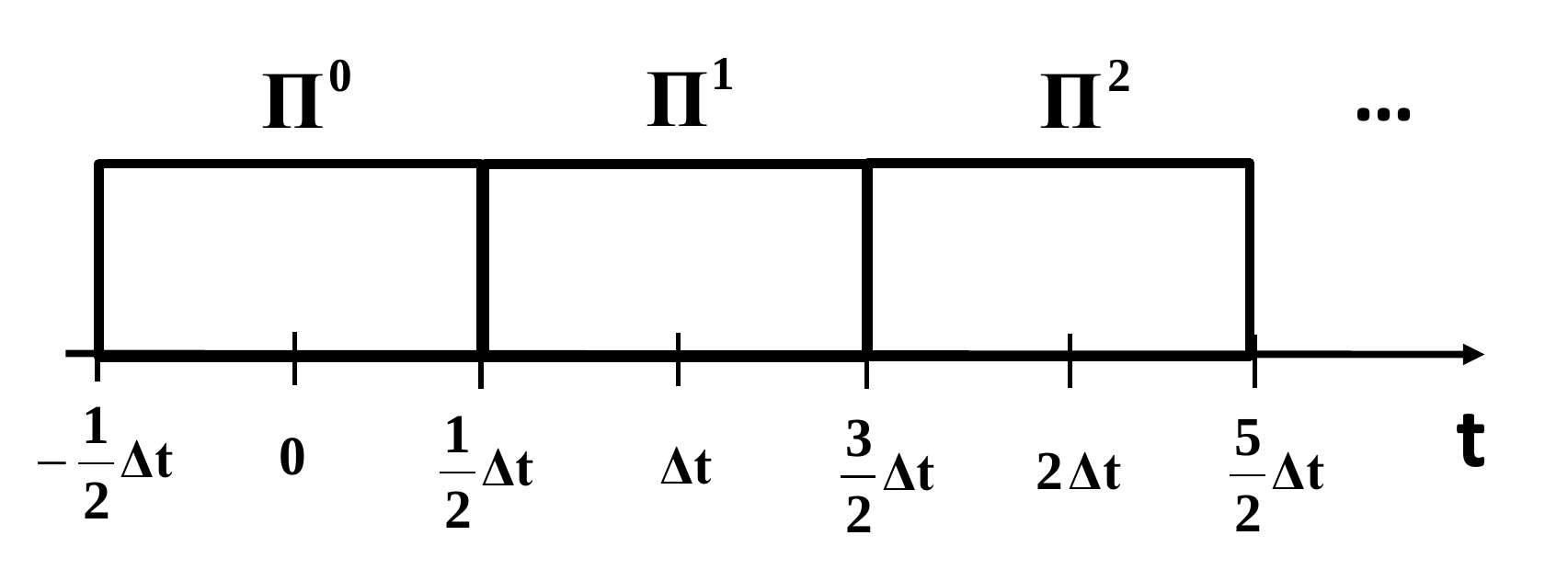}
    }\vspace{0.1 cm}
    \subfloat[\label{rooftop}]{%
      \includegraphics[width=2.9in]{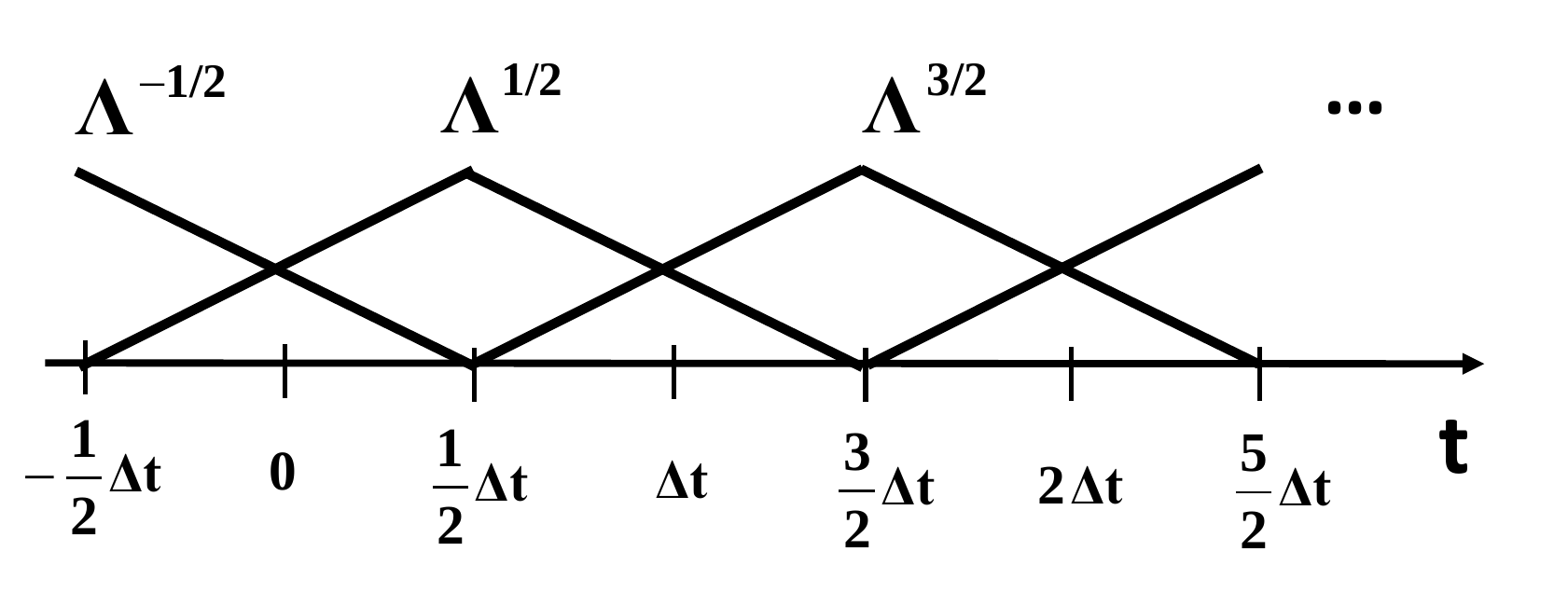}
    }
    \caption{Temporal basis functions of the lowest orders: (a) Piecewise constant (pulse) function and (b) Piecewise linear (rooftop) function.}
    \label{basis.time}
\end{figure}

%%%%%%%%%%%%%%%%%%%%%%%%%%%%%%%%%%%%%%%%%%%%%%%%%%%%%%%%%%%%%%%%%%%%%%%%%%%%%%%%%%%%%%%%%%%%%%%%%%%%%%%%%%%%%%
\subsection{Gather step}
%%%%%%%%%%%%%%%%%%%%%%%%%%%%%%%%%%%%%%%%%%%%%%%%%%%%%%%%%%%%%%%%%%%%%%%%%%%%%%%%%%%%%%%%%%%%%%%%%%%%%%%%%%%%%%
In the gather step, field values are interpolated at the positions of particles. Since Whitney basis functions are used to represent the field values, \eqref{E} and \eqref{B} can be directly used for the interpolation. Using \eqref{e.time.exp} and \eqref{b.time.exp}, $\mathbf{E}$ and $\mathbf{B}$ in their respective discrete times are expressed as
\begin{flalign}
&\mathbf{E}(\mathbf{r}_p,n\Delta t)
=\mathbf{E}^n(\mathbf{r}_p)
	= \sum^{N_e}_{i=1} e_i^n \mathbf{W}_i^1(\mathbf{r}_p), \\
&\mathbf{B}\left(\mathbf{r}_p,(n+1/2)\Delta t\right)
=\mathbf{B}^{n+\frac{1}{2}}(\mathbf{r}_p)
	= \sum^{N_f}_{i=1} b_i^{n+\frac{1}{2}} \mathbf{W}_i^2(\mathbf{r}_p),
\end{flalign}
where $\mathbf{r}_p$ is the position of the $p$-th particle.

%%%%%%%%%%%%%%%%%%%%%%%%%%%%%%%%%%%%%%%%%%%%%%%%%%%%%%%%%%%%%%%%%%%%%%%%%%%%%%%%%%%%%%%%%%%%%%%%%%%%%%%%%%%%%%
\subsection{Particle update}
%%%%%%%%%%%%%%%%%%%%%%%%%%%%%%%%%%%%%%%%%%%%%%%%%%%%%%%%%%%%%%%%%%%%%%%%%%%%%%%%%%%%%%%%%%%%%%%%%%%%%%%%%%%%%%
The next step is to update particle attributes such as position $\mathbf{r}_p(t)$ and velocity $\mathbf{v}_p(t)$. The equation of motion and Lorentz-Newton equation are utilized. For simplicity, we consider here a non-relativistic case:
\begin{flalign}
&\frac{d\mathbf{r}_p}{dt} = \mathbf{v}_p, \label{eq.motion}\\
&\frac{d\mathbf{v}_p}{dt} =
	\frac{q}{m}\left(\mathbf{E}+\mathbf{v}_p\times\mathbf{B}\right). \label{Lorentz.Newton}
\end{flalign}
In \eqref{Lorentz.Newton}, $q$ and $m$ are the charge and mass of the particle, respectively.
Using the leap-frog time update, \eqref{eq.motion} and \eqref{Lorentz.Newton} are discretized as
\begin{flalign}
\mathbf{r}_p^{n+1}-\mathbf{r}_p^n
	&=\Delta t \, \mathbf{v}_p^{n+\frac{1}{2}}, \label{eq.motion.dis}\\
\mathbf{v}_p^{n+\frac{1}{2}}-\mathbf{v}_p^{n-\frac{1}{2}}
	&=\frac{q\Delta t}{m}\left(\mathbf{E}^n + \mathbf{v}_p^n\times\mathbf{B}^n\right). \label{Lorentz.Newton.dis}
\end{flalign}
Note that \eqref{eq.motion.dis} is the explicit update whereas \eqref{Lorentz.Newton.dis} is implicit. Furthermore, the interpolated values of $\mathbf{v}_p^n$ and $\mathbf{B}^n$ in the right hand side of \eqref{Lorentz.Newton.dis} should be calculated a priori. We assume that $\mathbf{v}_p^n$ is expanded in the piecewise linear function because $\mathbf{v}_p^n$ and $\mathbf{B}^n$ are centered on the same time. Therefore, when $t=n\Delta t$,
\begin{flalign}
\mathbf{v}_p^n &= \frac{1}{2}\left(\mathbf{v}_p^{n+\frac{1}{2}} + \mathbf{v}_p^{n-\frac{1}{2}}\right), \\
\mathbf{B}^n &= \frac{1}{2}\left(\mathbf{B}^{n+\frac{1}{2}} + \mathbf{B}^{n-\frac{1}{2}}\right).
\end{flalign}
In contrast to $\mathbf{v}_p$, $\mathbf{B}$ is the function of space as well, so it is assumed that $\mathbf{B}^{n+\frac{1}{2}}(\mathbf{r}_p)=\mathbf{B}^{n+\frac{1}{2}}(\mathbf{r}_p^n)$.
Therefore, \eqref{Lorentz.Newton.dis} is modified to
\begin{flalign}
\mathbf{v}_p^{n+\frac{1}{2}}-\mathbf{v}_p^{n-\frac{1}{2}}&=
\frac{q\Delta t}{m}
\left[\mathbf{E}^n +
	\frac{1}{4}\mathbf{v}_p^{n+\frac{1}{2}}\times
		\left(\mathbf{B}^{n+\frac{1}{2}} + \mathbf{B}^{n-\frac{1}{2}}\right)
	+ \frac{1}{4}\mathbf{v}_p^{n-\frac{1}{2}}\times
		\left(\mathbf{B}^{n+\frac{1}{2}} + \mathbf{B}^{n-\frac{1}{2}}\right)
\right]. \label{Lorentz.Newton.m}
\end{flalign}
After some algebra, \eqref{Lorentz.Newton.m} can be succinctly expressed as
\begin{flalign}
\mathbf{v}_p^{n+\frac{1}{2}}
=\mathbf{N}^{-1} \cdot \mathbf{N}^T \cdot \mathbf{v}_p^{n-\frac{1}{2}}
	+ \frac{q\Delta t}{m}\mathbf{N}^{-1}\cdot\mathbf{E}^n, \label{Lorentz.Newton.implicit.succinct}
\end{flalign}
where
\begin{flalign}
\renewcommand{\arraystretch}{1.4}
\mathbf{N}=
    \begin{bmatrix}
                             1 & -\frac{q\Delta t}{2m}B_z^n &  \frac{q\Delta t}{2m}B_y^n \\
     \frac{q\Delta t}{2m}B_z^n &                          1 & -\frac{q\Delta t}{2m}B_x^n \\
    -\frac{q\Delta t}{2m}B_y^n &  \frac{q\Delta t}{2m}B_x^n &                          1
    \end{bmatrix},
\end{flalign}
\begin{flalign}
B_s^n=\frac{1}{2}\left(B_s^{n+\frac{1}{2}}+B_s^{n-\frac{1}{2}}\right), \quad s=x,\;y,\text{ or } z.
\end{flalign}
Note that $\mathbf{N}$ is unitless. In summary, \eqref{eq.motion.dis} together with \eqref{Lorentz.Newton.implicit.succinct} constitute the well-known (non-relativistic) equation of motion for a charged particle in an electromagnetic field. Note that the particle velocity should be updated before the update of the particle position.

%%%%%%%%%%%%%%%%%%%%%%%%%%%%%%%%%%%%%%%%%%%%%%%%%%%%%%%%%%%%%%%%%%%%%%%%%%%%%%%%%%%%%%%%%%%%%%%%%%%%%%%%%%%%%%
\subsection{Scatter step}
\label{scatter}
%%%%%%%%%%%%%%%%%%%%%%%%%%%%%%%%%%%%%%%%%%%%%%%%%%%%%%%%%%%%%%%%%%%%%%%%%%%%%%%%%%%%%%%%%%%%%%%%%%%%%%%%%%%%%%
This step is to assign charge density and current density back to the grid using the updated particle attributes. Of course, the fundamental question here is how to assign the particle charge to grid vertices (nodes) consistent to the assignment of the respective currents to grid edges. To achieve this, we use the same family of interpolatory functions for both, viz. Whitney functions. The advantage of using such functions is two-fold: (i) they preserve the total amount of the vertex-distributed charge for each particle and (ii) they exactly preserve the charge continuity equation (more on this below) and guarantee that no spurious charges arise in Gauss' law during the particle motion on the grid.

\begin{figure}[t]
	\centering
	\subfloat[\label{charge_scatter}]{%
      \includegraphics[width=2.0in]{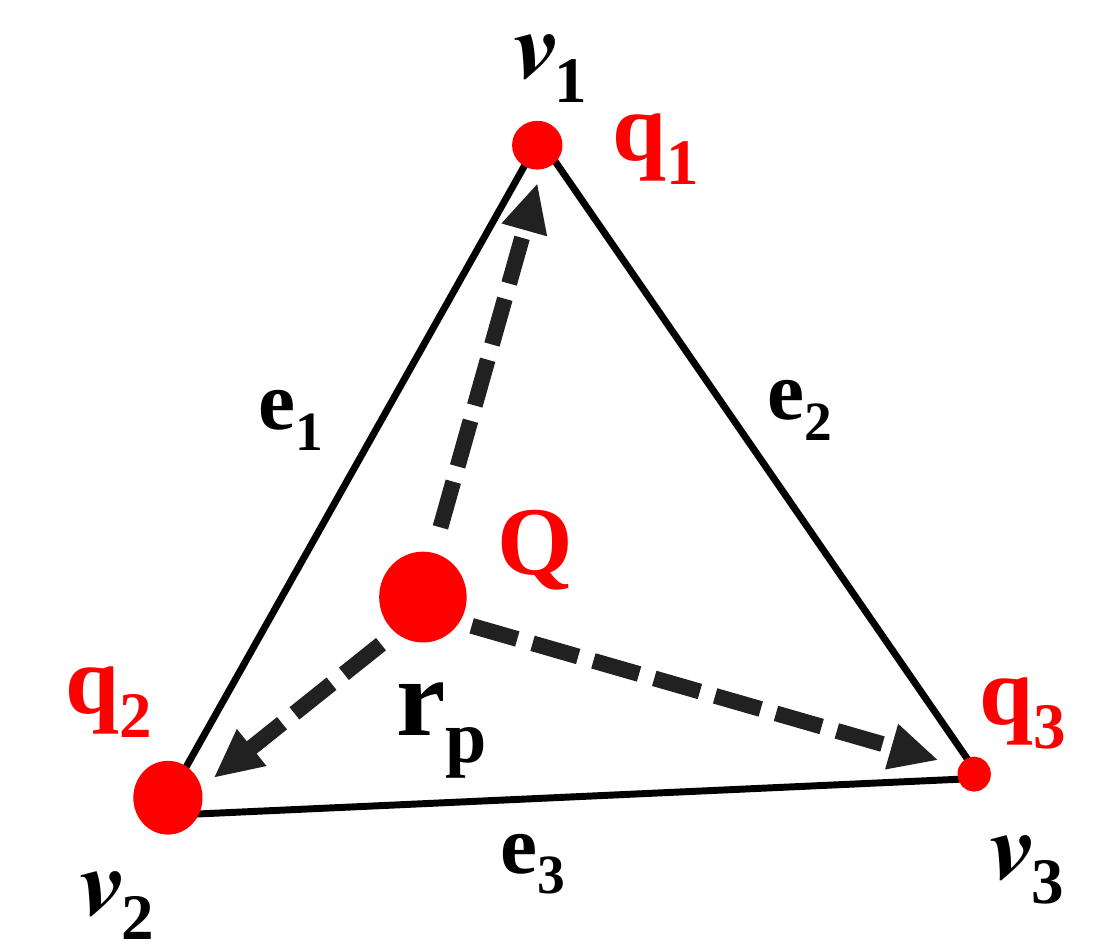}
    }
    \hspace{0.5 cm}
    \subfloat[\label{current_scatter}]{%
      \includegraphics[width=2.0in]{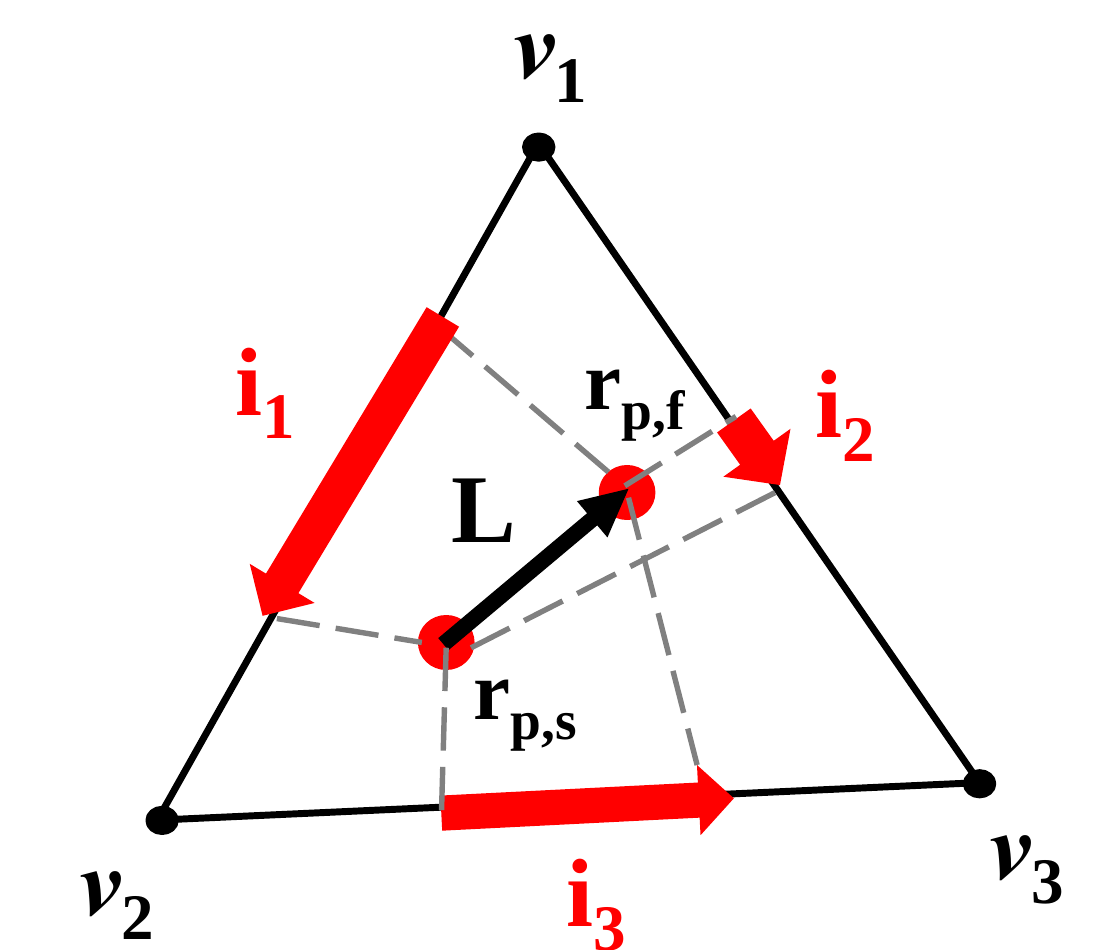}
    }
    \caption{Scatter step: (a) Nodal charge assignment from a charged particle placed at $\mathbf{r}_p$ with local numbering of vertices and edges and (b) Current assignment due to charge movement from $\mathbf{r}_{p,s}$ to $\mathbf{r}_{p,f}$ during $\Delta t$ with default directions for currents.}
    \label{Scatter.dis}
\end{figure}

To examine this, let us first consider the charge assignment. The charge $Q$ of the $p$-th particle is represented as a 0-form and distributed to the grid vertices so that
\begin{flalign}
q_i=Q W_i^0(\mathbf{r}_p)=Q\lambda_i(\mathbf{r}_p), \label{scatter.charge.assign}
\end{flalign}
where the subscript $i$ is the index of vertices and $W_i^0$ is the Whitney 0-form associated with vertex $i$. The value of the function $W_i^0(\mathbf{r}_p)$ is simply equal to the barycentric coordinate of the point $\mathbf{r}_p$ referred to the node $i$, i.e., $W_i^0(\mathbf{r}_p)=\lambda_i(\mathbf{r}_p)$ (see \ref{app.a}). When \eqref{scatter.charge.assign} is summed over all possible $i$ values,
\begin{flalign}
\sum_i q_i=\sum_i Q\lambda_i(\mathbf{r}_p)=Q\sum_i\lambda_i(\mathbf{r}_p)=Q,
\end{flalign}
 since $\sum_i \lambda_i(\mathbf{r}_p)=1$ holds (a partition of unity). The charge assignment \eqref{scatter.charge.assign} is illustrated in Fig. \ref{charge_scatter} with the local numbering of vertices and edges. Vertices and edges are represented by $\nu$ and $e$, respectively. Note that the charge values are only associated to the vertices of the triangle on which the particle is located.

Fig. \ref{current_scatter} describes the current assignment. The $p$-th particle of charge $Q$ moves from $\mathbf{r}_{p,s}$ to $\mathbf{r}_{p,f}$ during $\Delta t$ along straight path $\mathbf{L}$. For example, the current assigned to $e_1$ (edge 1) is
\begin{flalign}
i_1
=\frac{Q}{\Delta t} \int_{\mathbf{r}_{p,s}}^{\mathbf{r}_{p,f}} \mathbf{W}_1^1 (\mathbf{r}_p) \cdot d{\mathbf{L}}
=\frac{Q}{\Delta t} \left(\lambda_1^s \lambda_2^f - \lambda_1^f \lambda_2^s\right), \label{I.assign}
\end{flalign}
where $\lambda_i^s$ and $\lambda_i^f$ are shorthands of $\lambda_i(\mathbf{r}_{p,s})$ and $\lambda_i(\mathbf{r}_{p,f})$, respectively. See \ref{app.b} for further details on the evaluation of this line integral. We note that Equation (\ref{I.assign}) is a simpler version of the first equation in Section 3.2 of~\cite{Pinto14:Charge}, here evaluated in closed-form along a linear particle trajectory. The current values $i_2$ and $i_3$ can be obtained similarly.
During the scatter step $\Delta t$, particles might travel beyond a single triangle element and cross element edges.
In this case, the path can be simply divided into smaller segments whereby each segment resides entirely within a single triangle. The scatter step above can then be applied to each segment.

%%%%%%%%%%%%%%%%%%%%%%%%%%%%%%%%%%%%%%%%%%%%%%%%%%%%%%%%%%%%%%%%%%%%%%%%%%%%%%%%%%%%%%%%%%%%%%%%%%%%%%%%%%%%%%
\subsection{Charge conservation}
\label{verification}
%%%%%%%%%%%%%%%%%%%%%%%%%%%%%%%%%%%%%%%%%%%%%%%%%%%%%%%%%%%%%%%%%%%%%%%%%%%%%%%%%%%%%%%%%%%%%%%%%%%%%%%%%%%%%%
To verify charge conservation, let us consider the semi-discrete continuity equation
\begin{flalign}
\frac{d}{dt}\mathbf{q}+\widetilde{\mathbf{S}} \cdot \mathbf{i}=0, \label{semi.dis.cont}
\end{flalign}
where the array $\mathbf{q}$ represents the amount of charge at all vertices, i.e.,
$\mathbf{q}=[q_1(t),q_2(t),\cdots ,q_{N_\nu}(t)]^T$,
$N_\nu$ being the number of vertices in the grid, and $\widetilde{\mathbf{S}}$ being the incidence matrix
associated with the (discrete) divergence operator in the dual grid
~\cite{He06:Geometric, Teixeira99:Lattice, Clemens01:Discrete, He05:Degrees}.
Note that, similarly to $\mathbf{C}$, all elements of $\widetilde{\mathbf{S}}$ are in the set \{-1,0,1\}. Applying a leap-frog time update to \eqref{semi.dis.cont}, we obtain
\begin{flalign}\label{full.dis.cont}
\frac{\mathbf{q}^{n+1}-\mathbf{q}^n}{\Delta t}+\widetilde{\mathbf{S}} \cdot \mathbf{i}^{n+\frac{1}{2}}=0.
\end{flalign}
Then, let us consider $\nu_1$ (vertex 1) without loss of generality. The time rate of charge variation at $\nu_1$ is
\begin{flalign}\label{Q.node1}
\frac{q_1^{n+1} - q_1^n}{\Delta t}=\frac{Q\lambda_1^f}{\Delta t} - \frac{Q\lambda_1^s}{\Delta t}
=\frac{Q}{\Delta t}(\lambda_1^f-\lambda_1^s).
\end{flalign}
On the other hand, the current flowing out of $\nu_1$ can be computed as
\begin{flalign}
(\widetilde{\mathbf{S}}\mathbf{i}^{n+\frac{1}{2}})_1&=i_1+i_2 \notag\\
&=\frac{Q}{\Delta t}
	\left[
		\int_{\mathbf{r}_{p,s}}^{\mathbf{r}_{p,f}} \mathbf{W}_1^1 (\mathbf{r}_p) \cdot d{\mathbf{L}}
		+\int_{\mathbf{r}_{p,s}}^{\mathbf{r}_{p,f}} \mathbf{W}_2^1 (\mathbf{r}_p) \cdot d{\mathbf{L}}
	\right] \notag\\
&=\frac{Q}{\Delta t}
	\left[
		\left(\lambda_1^s \lambda_2^f - \lambda_1^f \lambda_2^s\right)
		+\left(\lambda_1^s \lambda_3^f - \lambda_1^f \lambda_3^s\right)
	\right] =\frac{Q}{\Delta t}
	\left[\lambda_1^s - \lambda_1^f \right], \label{I.node1}
\end{flalign}
where the property $\lambda_1+\lambda_2+\lambda_3=1$ has been used, and the Whitney edge basis functions are indexed in an ascending order fashion (instead of a cyclic order) such that
\begin{flalign}
\mathbf{W}_1^1 (\mathbf{r}_p) &= \lambda_1(\mathbf{r}_p)\na\lambda_2(\mathbf{r}_p) - \lambda_2(\mathbf{r}_p)\na\lambda_1(\mathbf{r}_p),  \\
\mathbf{W}_2^1 (\mathbf{r}_p) &= \lambda_1(\mathbf{r}_p)\na\lambda_3(\mathbf{r}_p) - \lambda_3(\mathbf{r}_p)\na\lambda_1(\mathbf{r}_p),  \\
\mathbf{W}_3^1 (\mathbf{r}_p) &= \lambda_2(\mathbf{r}_p)\na\lambda_3(\mathbf{r}_p) - \lambda_3(\mathbf{r}_p)\na\lambda_2(\mathbf{r}_p).
\end{flalign}
As the sum of \eqref{Q.node1} and \eqref{I.node1} equals zero, the continuity equation is verified exactly.

\begin{figure}[t]
	\centering
	\subfloat[\label{Whitney_0form}]{%
      \includegraphics[width=2.0in]{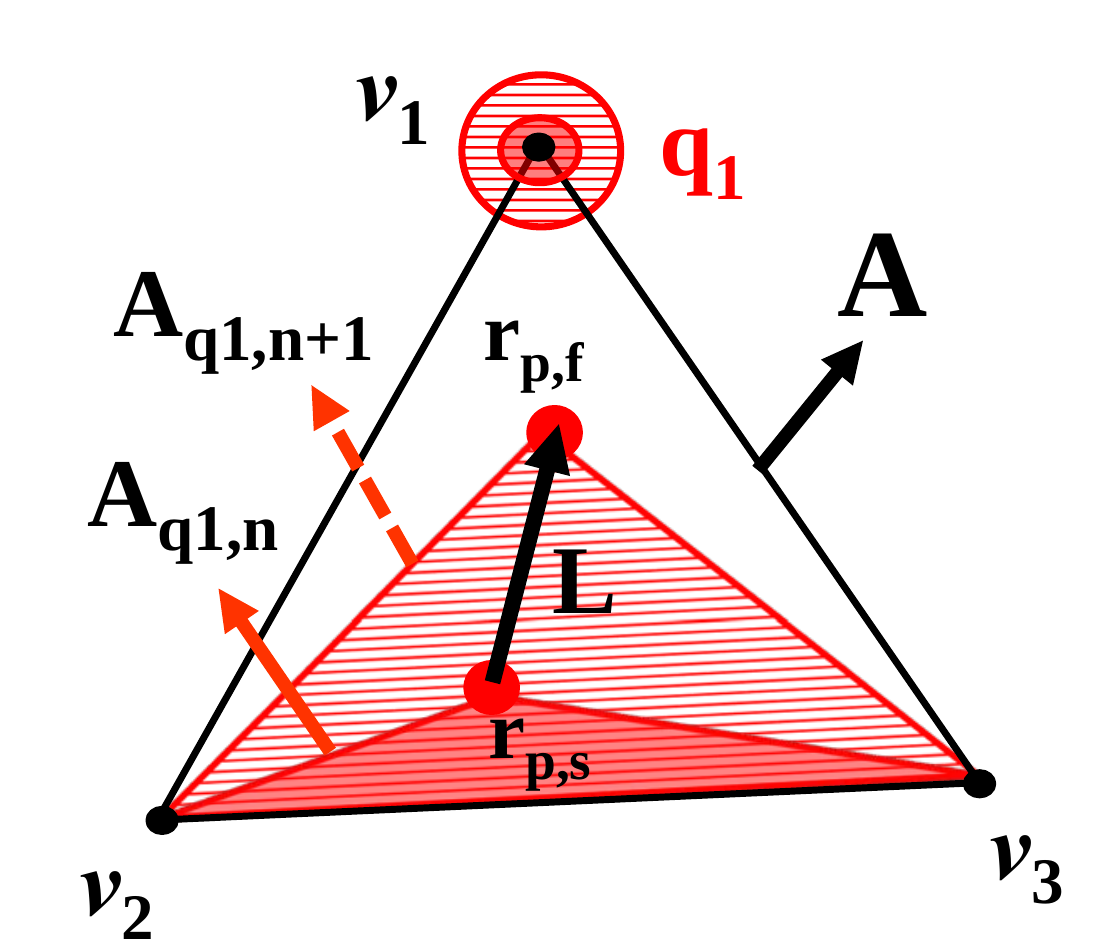}
    }
    \hspace{0.5 cm}
    \subfloat[\label{Whitney_1form}]{%
      \includegraphics[width=2.0in]{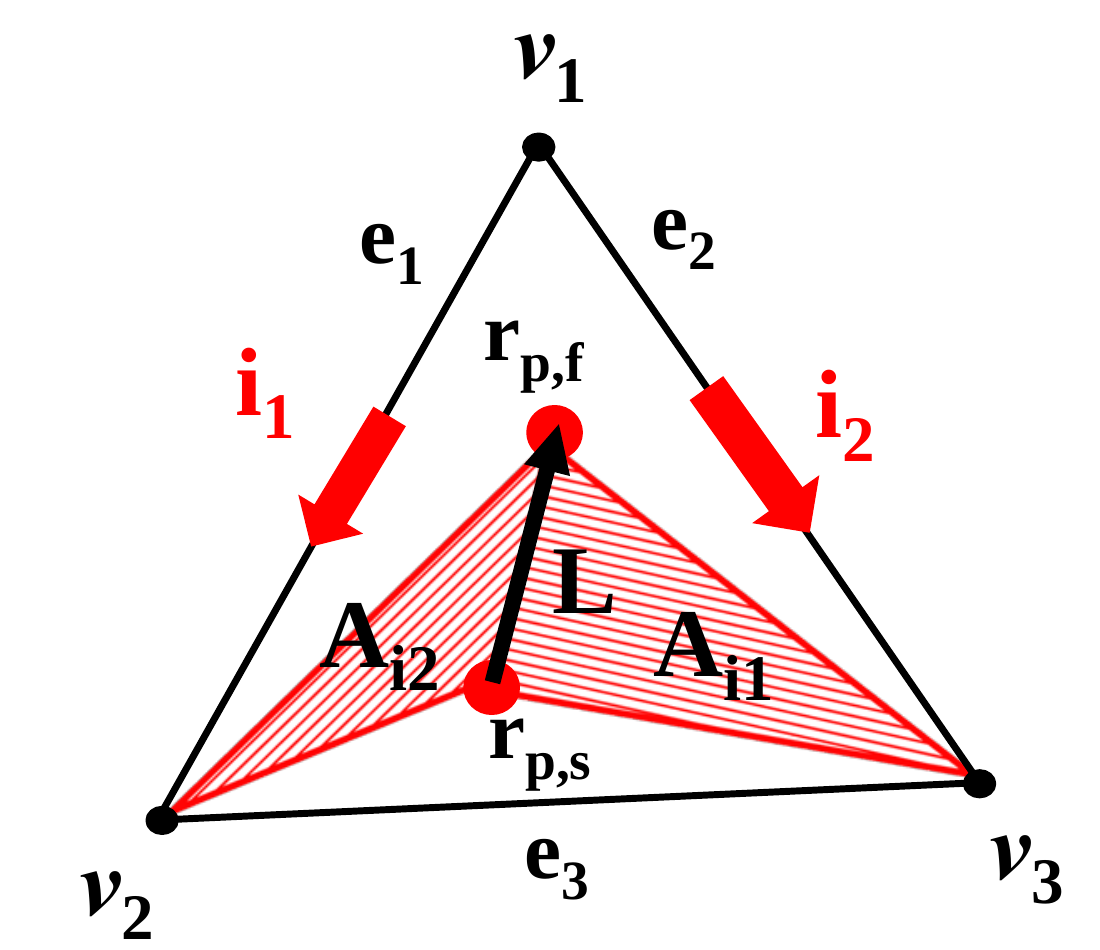}
    }
    \caption{Geometric representation of charge-conservation identity \eqref{A.area}: (a) Variation of Whitney 0-forms coefficients (barycentric coordinates) associated with $\nu_1$ (vertex 1) during a time interval $\Delta t$, where the solid red region indicates the starting time instant and the striped red region (which includes the solid red region) indicates the finishing time instant. (b) Areas associated with the magnitude of induced currents (as computed by Whitney 1-forms) on adjacent edges $e_1$ (edge 1) and $e_2$ (edge 2) during $\Delta t$
 (see also 
    Figs.~\ref{bary.coord} and~\ref{Whitney.1form.area}).
It is clear that $A_{q1,n+1} - A_{q1,n} = A_{i1} + A_{i2}$.
}
    \label{Whitney.forms}
\end{figure}

% begin red-colored text
The above derivation can be interpreted geometrically by understanding the geometric representation of Whitney 0-forms and 1-forms. Let us consider $\nu_1$ again. As explained in \ref{app.a}
and illustrated in Fig.~\ref{bary.coord}, barycentric coordinates can be visualized as a ratio of two areas.
The variation on the charge assigned to $\nu_1$ during $\Delta t$ is illustrated in terms of such areas in Fig.~\ref{Whitney_0form} and is expressed as
\begin{flalign}
\frac{Q}{\Delta t}(\lambda_1^f-\lambda_1^s)=\frac{Q}{\Delta t}\frac{A_{q1,n+1}-A_{q1,n}}{A}, \label{Whitney.0form.q}
\end{flalign}
where $A_{q1,n1}$ and $A_{q1,n+1}$ are the triangle areas as indicated in Fig.~\ref{Whitney_0form}, and $A$ is the area of the whole triangle (grid element) defined by $\nu_1$, $\nu_2$ and $\nu_3$. On the other hand, the current flowing out of $\nu_1$ is the sum of the currents along the edges touching $\nu_1$, that is $e_1$ and $e_2$. Referring to the geometric interpretation of the integral of Whitney 1-forms provided in \ref{app.a}, the sum of these two currents is evaluated as
\begin{flalign}
i_1+i_2
=\frac{Q}{\Delta t}
	\left[
		-\int_{\mathbf{r}_{p,s}}^{\mathbf{r}_{p,f}} \mathbf{W}_1^1 (\mathbf{r}_p) \cdot d{\mathbf{L}}
		-\int_{\mathbf{r}_{p,s}}^{\mathbf{r}_{p,f}} \mathbf{W}_2^1 (\mathbf{r}_p) \cdot d{\mathbf{L}}
	\right]
=-\frac{Q}{\Delta t}
	\left[\frac{A_{i1}}{A} + \frac{A_{i2}}{A}\right], \label{Whitney.1form.I}
\end{flalign}
where $A_{i1}$ and $A_{i2}$ are indicated in Fig.~\ref{Whitney_1form}, and the minus sign is due to the relative orientations of the path $\mathbf{L}$ and Whitney 1-forms.  Since
\begin{flalign}
A_{q1,n+1} - A_{q1,n} = A_{i1} + A_{i2},
\label{A.area}
\end{flalign}
the sum of \eqref{Whitney.0form.q} and \eqref{Whitney.1form.I} is identically zero.

%%%%%%%%%%%%%%%%%%%%%%%%%%%%%%%%%%%%%%%%%%%%%%%%%%%%%%%%%%%%%%%%%%%%%%%%%%%%%%%%%%%%%%%%%%%%%%%%%%%%%%%%%%%%%%
\subsection{Gauss' law preservation}
\label{gauss}
%%%%%%%%%%%%%%%%%%%%%%%%%%%%%%%%%%%%%%%%%%%%%%%%%%%%%%%%%%%%%%%%%%%%%%%%%%%%%%%%%%%%%%%%%%%%%%%%%%%%%%%%%%%%%%
We next demonstrate that Gauss' law is automatically satisfied for all time steps if proper initial conditions are used. By left-multiplying  both sides of \eqref{e.update} by the discrete divergence matrix $\widetilde{\mathbf{S}}$, we obtain
\begin{flalign}
\widetilde{\mathbf{S}}\cdot\left[\star_{\epsilon}\right]\cdot
	\left(\frac{\mathbf{e}^{n+1} - \mathbf{e}^n}{\Delta t}\right)
	= \widetilde{\mathbf{S}}\cdot\mathbf{C}^T\cdot\left[\star_{\mu^{-1}}\right]\cdot\mathbf{b}^{n+\frac{1}{2}}
	-\widetilde{\mathbf{S}}\cdot\mathbf{i}^{n+\frac{1}{2}}. \label{div.D}
\end{flalign}
% begin red-colored text
The first term of the right-hand side of \eqref{div.D} vanishes due to the exact sequence property for the dual grid, i.e., $\widetilde{\mathbf{S}}\cdot\mathbf{C}^T=0$~\cite{Teixeira99:Lattice, Clemens01:Discrete, Flanders:Differential}\footnote{The identity $\widetilde{\mathbf{S}} \cdot \mathbf{C}^T=0$ can be recognized as the discrete analogue of $\na\cdot\na\times = 0$.}.
% end red-colored text
Using the discrete continuity equation \eqref{full.dis.cont}, we can rearrange \eqref{div.D} as
\begin{flalign}
\widetilde{\mathbf{S}}\cdot\left[\star_{\epsilon}\right]\cdot
	\left(\frac{\mathbf{e}^{n+1} - \mathbf{e}^n}{\Delta t}\right)
	= \frac{\mathbf{q}^{n+1}-\mathbf{q}^n}{\Delta t}, \label{div.D.modified}
\end{flalign}
which is the discrete version of
\begin{flalign}
\frac{\pa}{\pa t}\na\cdot \mathbf{D} = \frac{\pa}{\pa t}\rho. \label{div.D.cont}
\end{flalign}
Therefore, Gauss' law is preserved for all time steps if the initial condition
$\widetilde{\mathbf{S}}\cdot\left[\star_{\epsilon}\right]\cdot\mathbf{e}^0=\mathbf{q}^0$ is met.

% begin red-colored text
For completeness, we show next that Gauss' law for magnetism is also satisfied if appropriate initial conditions are used. 
By taking discrete divergence matrix $\mathbf{S}$ in both sides of \eqref{b.update}, we have\footnote{Note that $\mathbf{S}$ is distinct from $\widetilde{\mathbf{S}}$ since $\mathbf{S}$ refers to the primal grid (i.e., the computational mesh itself) whereas $\widetilde{\mathbf{S}}$ refers to the dual grid (See~\cite{Teixeira99:Lattice, Clemens01:Discrete}).}
\begin{flalign}
\mathbf{S}\cdot
	\left(\frac{\mathbf{b}^{n+\frac{1}{2}} - \mathbf{b}^{n-\frac{1}{2}}}{\Delta t}\right)
	= -\mathbf{S}\cdot\mathbf{C}\cdot\mathbf{e}^n=0, \label{div.D.mag}
\end{flalign}
where the second equality is from the exact sequence property in the primal grid, i.e., $\mathbf{S}\cdot\mathbf{C}=0$. The relation \eqref{div.D.mag} is the discrete version of
\begin{flalign}
\frac{\pa}{\pa t}\na\cdot \mathbf{B} = 0. \label{div.D.mag.cont}
\end{flalign}
Therefore, Gauss' law for magnetism is also preserved for all times if $\mathbf{b}^0$ is such that $\mathbf{S}\cdot\mathbf{b}^0=0$.
% end red-colored text

%%%%%%%%%%%%%%%%%%%%%%%%%%%%%%%%%%%%%%%%%%%%%%%%%%%%%%%%%%%%%%%%%%%%%%%%%%%%%%%%%%%%%%%%%%%%%%%%%%%%%%%%%%%%%%
\subsection{Time-update sequence}
\begin{figure}[t]
    \centering
    \includegraphics[width=3.6in]{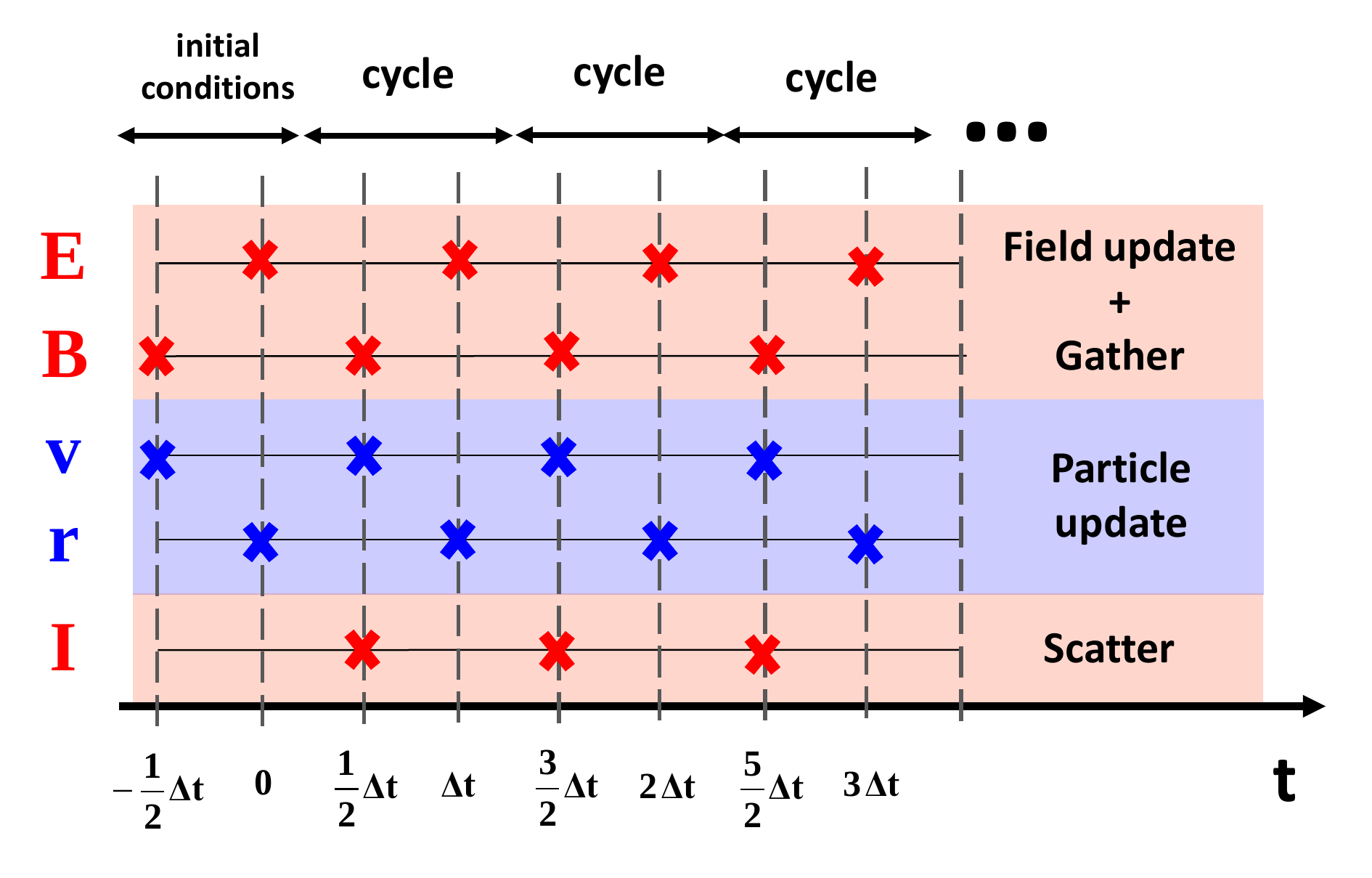}
    \caption{Full time-update procedure for the charge-conserving PIC algorithm.}
    \label{full.time.update}
\end{figure}

%%%%%%%%%%%%%%%%%%%%%%%%%%%%%%%%%%%%%%%%%%%%%%%%%%%%%%%%%%%%%%%%%%%%%%%%%%%%%%%%%%%%%%%%%%%%%%%%%%%%%%%%%%%%%%
Using the above equations, the overall time-update procedure is carried out in the following sequence. Initial conditions for $\mathbf{E}^0$, $\mathbf{B}^{-\frac{1}{2}}$, $\mathbf{v}_p^{-\frac{1}{2}}$, and $\mathbf{r}_p^0$ are first assumed. During each cycle, $\mathbf{b}^{n+\frac{1}{2}}$ is first calculated. Then, $\mathbf{E}^n$ and $\mathbf{B}^{n+\frac{1}{2}}$ are interpolated at particle positions. Next, after the particle acceleration $\mathbf{v}_p^{n+\frac{1}{2}}$ is performed, the particle push $\mathbf{r}_p^{n+1}$ is performed for all particles.  Next, currents $\mathbf{i}^{n+\frac{1}{2}}$ are assigned (scattered) to grid edges. Finally, $\mathbf{e}^{n+1}$ is updated. Note that $\mathbf{v}_p$ and $\mathbf{r}_p$ are 3$\times$1 column vectors. The procedure is illustrated in Fig. \ref{full.time.update} and each step is enumerated below.

\allowdisplaybreaks
% B update
\begin{flalign}
&\text{1) } \mathbf{B} \text{ update : }
\mathbf{b}^{n+\frac{1}{2}}
	=\mathbf{b}^{n-\frac{1}{2}} - \Delta t \, \mathbf{C}\cdot \mathbf{e}^n \notag\\
% Gather - E
&\text{2) } \mathbf{E} \text{ gather : }
\mathbf{E}^n = \sum^{N_e}_{i=1} e_i^n \mathbf{W}_i^1(\mathbf{r}_p^n) \notag\\
% Gather - B
&\text{3) } \mathbf{B} \text{ gather : }
\mathbf{B}^{n+\frac{1}{2}}	= \sum^{N_f}_{i=1} b_i^{n+\frac{1}{2}} \mathbf{W}_i^2(\mathbf{r}_p^n) \notag\\
% particle acceleration (particle velocity update)
&\text{4) } \text{Particle acceleration : }
\mathbf{v}_p^{n+\frac{1}{2}}=\mathbf{N}^{-1}\cdot\mathbf{N}^T\cdot\mathbf{v}_p^{n-\frac{1}{2}}
	+ \frac{q\Delta t}{m}\mathbf{N}^{-1}\cdot\mathbf{E}^n \notag\\
% particle push (particle position update)
&\text{5) } \text{Particle push : }
\mathbf{r}_p^{n+1} = \mathbf{r}_p^n + \Delta t \mathbf{v}_p^{n+\frac{1}{2}} \notag\\
% current scatter
&\text{6) } \mathbf{I} \text{ scatter : }
i_i^{n+\frac{1}{2}}
	=\frac{Q}{\Delta t} \int_{\mathbf{r}_{p,s}}^{\mathbf{r}_{p,f}} \mathbf{W}_i^1(\mathbf{r}_p)\cdot d{\mathbf{L}} \notag\\
% E update
&\text{7) } \mathbf{E} \text{ update : }
\left[\star_{\epsilon}\right] \cdot \mathbf{e}^{n+1}
	=\left[\star_{\epsilon}\right] \cdot \mathbf{e}^n + \Delta t
		\left(
		 \mathbf{C}^T \cdot \left[\star_{\mu^{-1}}\right]\cdot\mathbf{b}^{n+\frac{1}{2}}	-\mathbf{i}^{n+\frac{1}{2}}
		\right) \notag
\end{flalign}
The algorithm utilizes an ``intelligent'' mesh for tracking particles at each time step without resorting to iterative search or lookup tables. The intelligent mesh is constructed (initialized) once the input mesh is loaded and it includes (adds) the necessary connectivity information among mesh elements to efficiently determine the element location of each particle in the next time step. The computing time of this process is minimal because it always starts the particle search from adjacent elements.
It should also be stressed that the proposed scatter-gather algorithm is independent of the time integration scheme and it can be combined with other schemes as well.

%%%%%%%%%%%%%%%%%%%%%%%%%%%%%%%%%%%%%%%%%%%%%%%%%%%%%%%%%%%%%%%%%%%%%%%%%%%%%%%%%%%%%%%%%%%%%%%%%%%%%%%%%%%%%%
\section{Validation}
\label{schem.val}
\begin{figure}[t]
	\centering
	\subfloat[\label{case1_t0}]{%
      \includegraphics[width=1.8in]{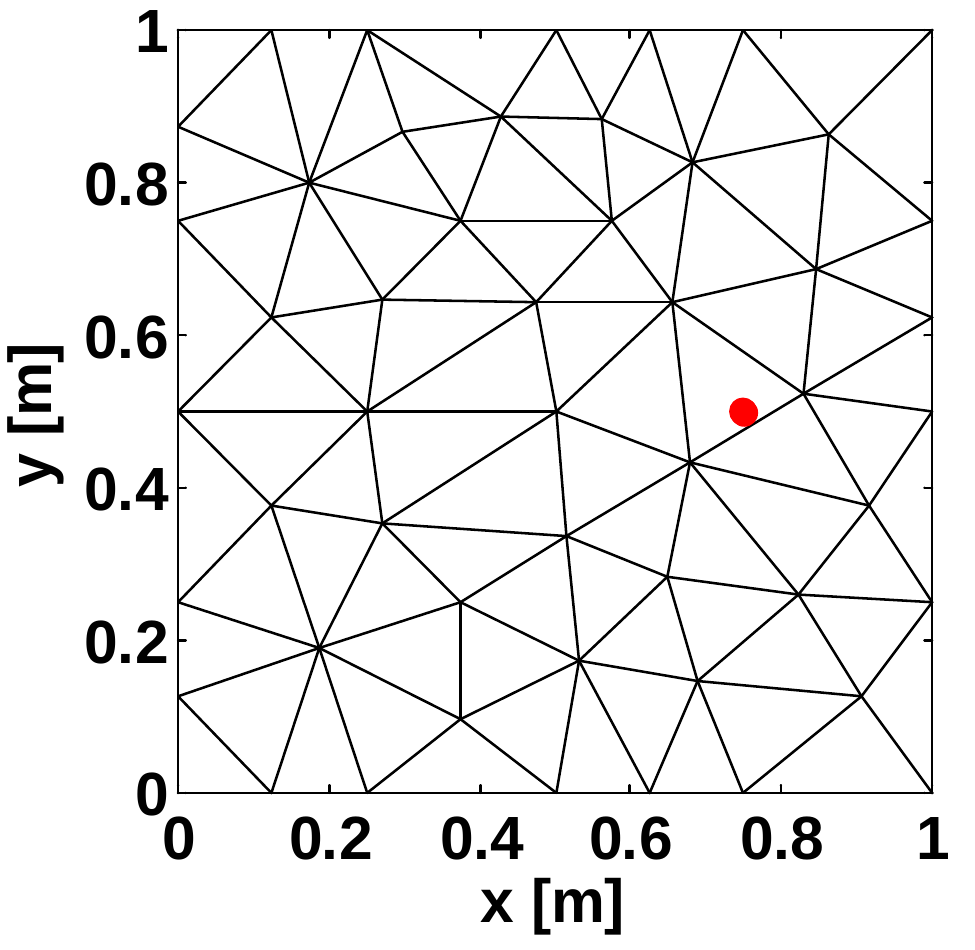}
    }
    \hspace{0.5 cm}
    \subfloat[\label{case1_t50}]{%
      \includegraphics[width=1.8in]{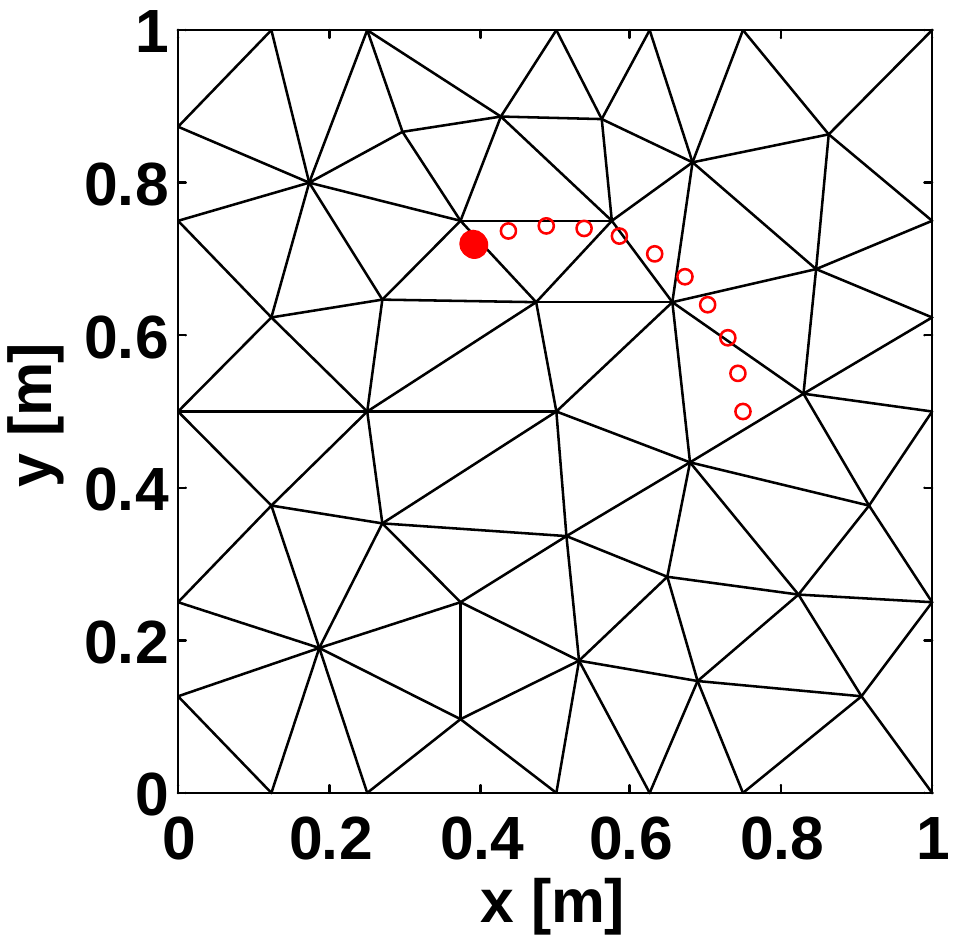}
    }
    \\
   	\subfloat[\label{case1_t100}]{%
      \includegraphics[width=1.8in]{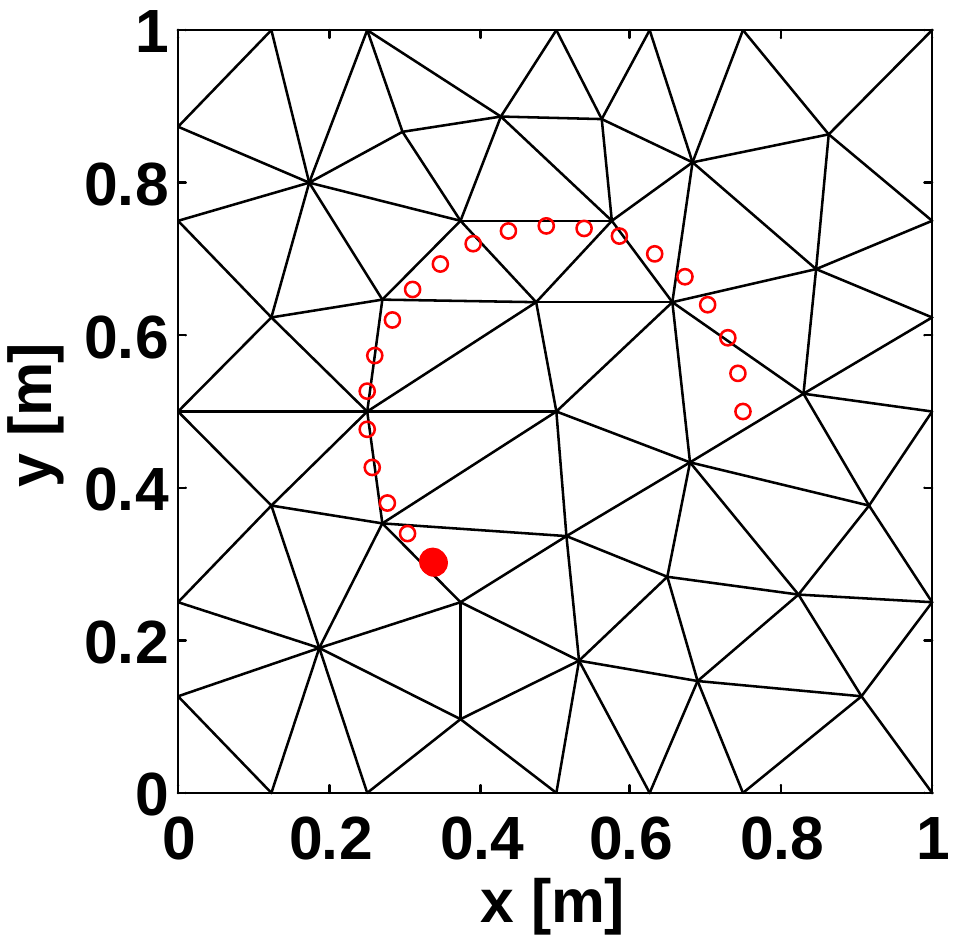}
    }
    \hspace{0.5 cm}
    \subfloat[\label{case1_t200}]{%
      \includegraphics[width=1.8in]{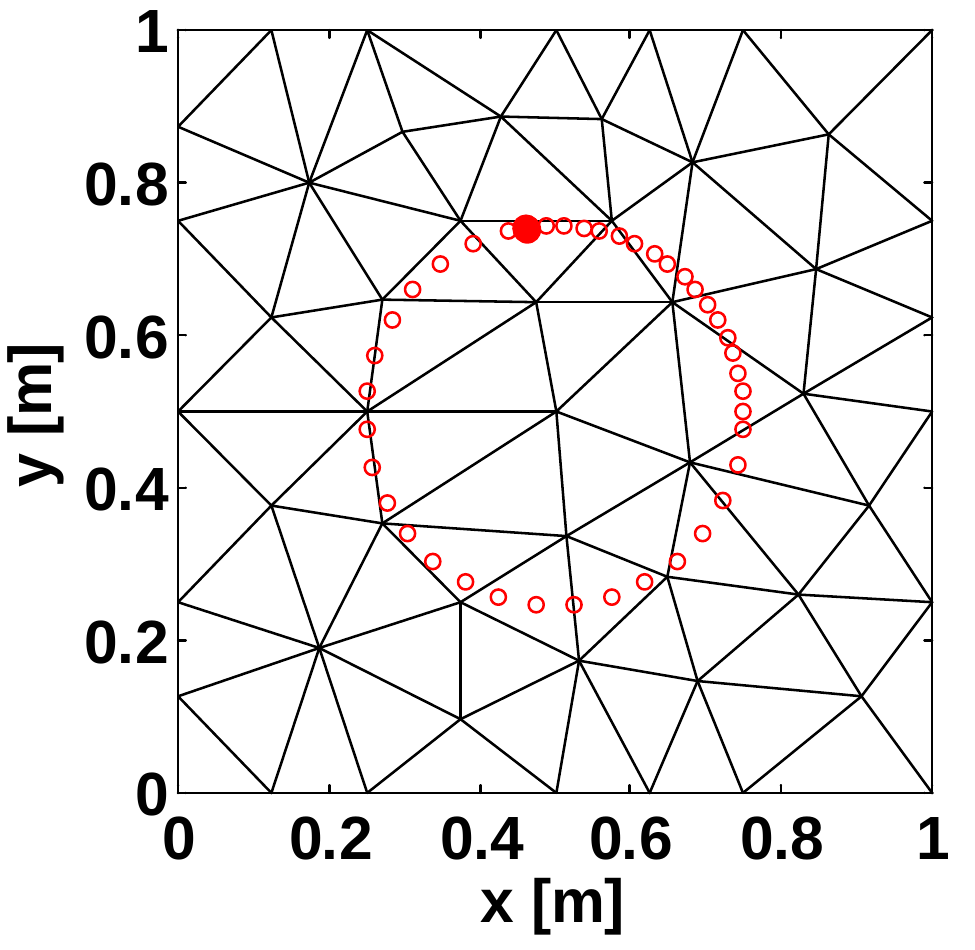}
    }
    \caption{Movement of a single particle in the uniform static magnetic field at different time instants ($\Delta t$ = 0.1 ns): (a) t = 0, (b) t = 50$\Delta t$, (c) t = 100$\Delta t$, and (d) t = 200$\Delta t$.}
    \label{particle.move}
\end{figure}

\begin{figure}[t]
	\centering
	\subfloat[\label{charge.total}]{%
      \includegraphics[width=2.6in]{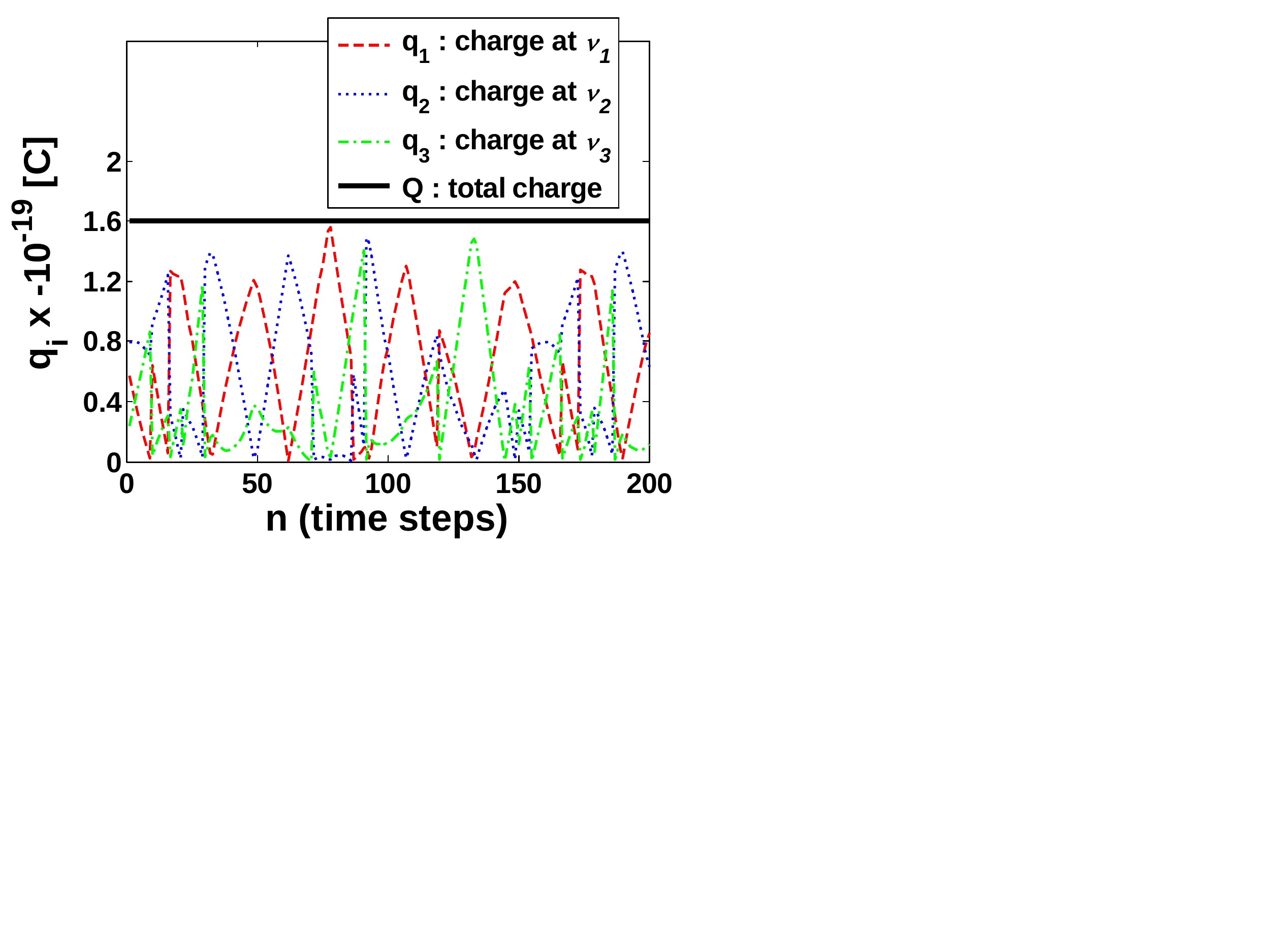}
    }
	\subfloat[\label{velocity}]{%
      \includegraphics[width=2.6in]{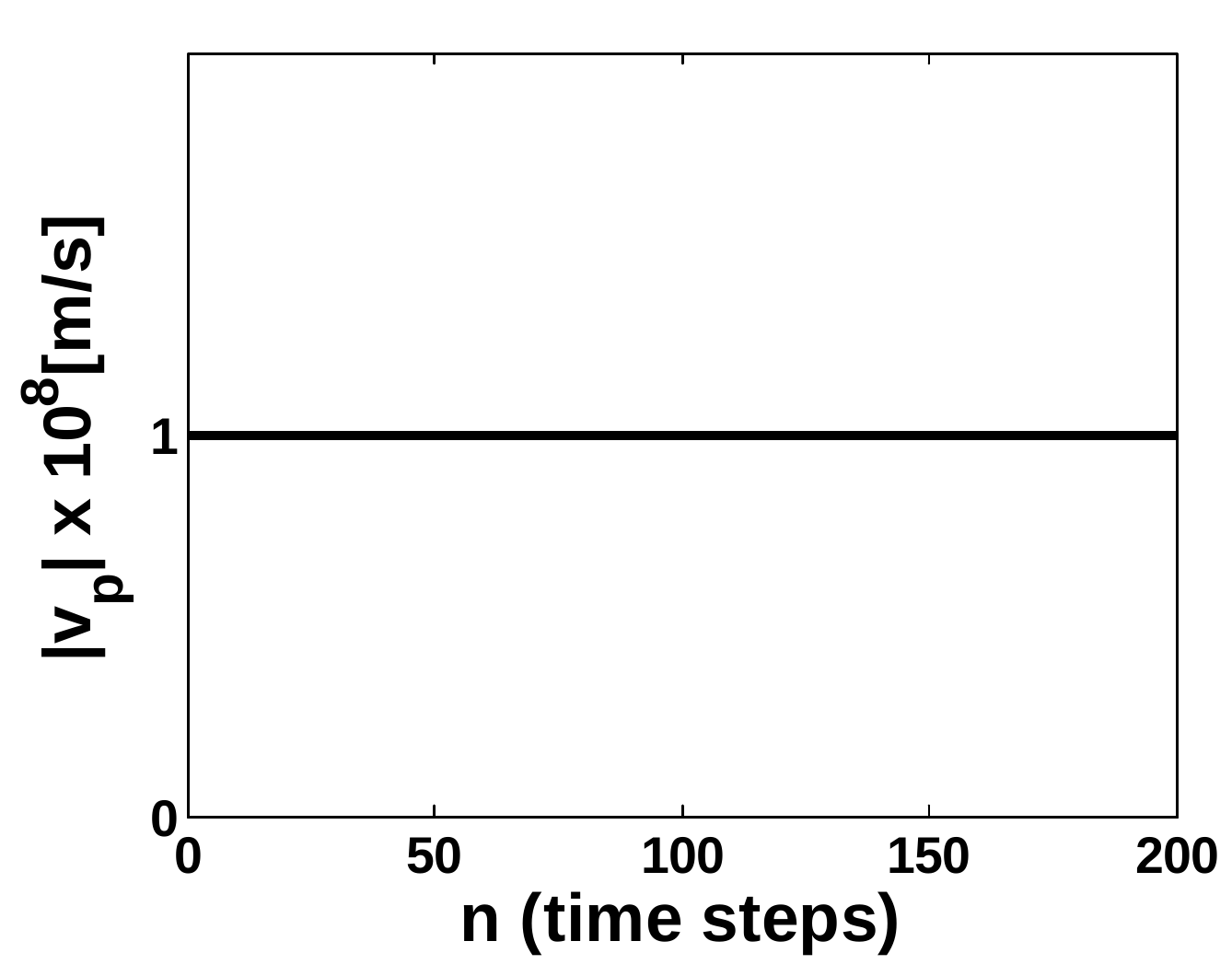}
    }
    \caption{Charge and energy conservation: (a) Distributed amounts of charge to local vertices and their sum at all time steps and (b) Absolute value of the particle velocity at all time steps.}
    \label{charge.energy.conservation}
\end{figure}

%%%%%%%%%%%%%%%%%%%%%%%%%%%%%%%%%%%%%%%%%%%%%%%%%%%%%%%%%%%%%%%%%%%%%%%%%%%%%%%%%%%%%%%%%%%%%%%%%%%%%%%%%%%%%%
% begin red-colored text
Let us consider a simple cyclotron motion for which a uniform static magnetic field is excited along the $z$-direction.
The static magnetic flux density is $B_z=2.275\times10^{-3}$ Wb/m$^{2}$, which produces the gyroradius of 0.25 m using $B_z=(mv)/(rq)$, where $m=9.1\times 10^{-31}$ kg, $v=10^8$ m/s, and $q=-1.6\times 10^{-19}$ C.
% end red-colored text
Fig. \ref{particle.move} shows the snapshots of the movement of a single particle at selected time steps. As the scheme is conditionally stable, time step should be less than the Courant limit $\Delta t_c=0.14887$ ns, which is the function of the mesh element sizes and is computed from the maximum eigenvalue of the stiffness matrix~\cite{Moon14:Trade}. This $\Delta t_c$ is less than $\Delta l/|\mathbf{v_p}|\approx 0.1/10^{8} = 10^{-9}$ s, where $\Delta l$ is the typical edge length of the triangular grid elements. It can be observed that the particle exactly shows the circular motion of 0.25 m radius. In this simulation, a pair of particles with the opposite charges are initially placed in the same location, so that net charge density and electric fields are initially zero. In contrast to the negatively charged particle, the particle with positive charge is assumed to be stationary due to its much larger mass, which is not shown in Fig. \ref{particle.move}.
% begin red-colored text
Fig. \ref{charge.total} and \ref{velocity} show the amount of vertex-distributed charge and the absolute value of the particle velocity as a function of time, respectively. The total charge remains constant by the virtue of the consistent particle interpolation in the scatter-gather algorithm. The absolute value of the particle velocity (hence, energy) also remains constant as well due to a negligible electric field. Table \ref{T.charge.cons} shows similar results as Fig. \ref{charge.energy.conservation}, but extended up to $10^6$ time steps to further verify charge and energy conservation.
% end red-colored text

\begin{table}[t]
\begin{center}
\renewcommand{\arraystretch}{1.3}
\caption{Charge and energy conservation at large time steps.}
\begin{tabular}{cccccc}
    \hline
    $n$ & $q_1$ & $q_2$ & $q_3$ & $Q$ & $|v_p|$ \\
    \hline
    $10^1$ & -6.410056 $\times 10^{-20}$ & -9.245821 $\times 10^{-20}$ & -3.441224 $\times 10^{-21}$ &
		-1.600000 $\times 10^{-19}$ & 9.999999 $\times 10^{7}$ \\
	$10^2$ & -7.635486 $\times 10^{-20}$ & -7.154041 $\times 10^{-20}$ & -1.210471 $\times 10^{-20}$ &
	    -1.600000 $\times 10^{-19}$ & 9.999999 $\times 10^{7}$ \\
    $10^3$ & -6.187120 $\times 10^{-20}$ & -7.721123 $\times 10^{-20}$ & -2.091755 $\times 10^{-20}$ &
		-1.600000 $\times 10^{-19}$ & 9.999999 $\times 10^{7}$ \\
    $10^4$ & -5.772639 $\times 10^{-21}$ & -1.472014 $\times 10^{-19}$ & -7.025898 $\times 10^{-21}$ &
		-1.600000 $\times 10^{-19}$ & 9.999999 $\times 10^{7}$ \\
    $10^5$ & -5.766949 $\times 10^{-20}$ & -2.809120 $\times 10^{-20}$ & -7.423930 $\times 10^{-20}$ &
		-1.600000 $\times 10^{-19}$ &  1.000000 $\times 10^{8}$ \\	
	$10^6$ & -1.480969 $\times 10^{-20}$ & -1.365091 $\times 10^{-21}$ & -1.438252 $\times 10^{-19}$ &
		-1.600000 $\times 10^{-19}$ &  1.000000 $\times 10^{8}$ \\	
	\hline
\end{tabular}
\label{T.charge.cons}
\end{center}
\end{table}					
					
\begin{figure}[t]
	\centering
	\subfloat[\label{case2_t0}]{%
      \includegraphics[width=1.8in]{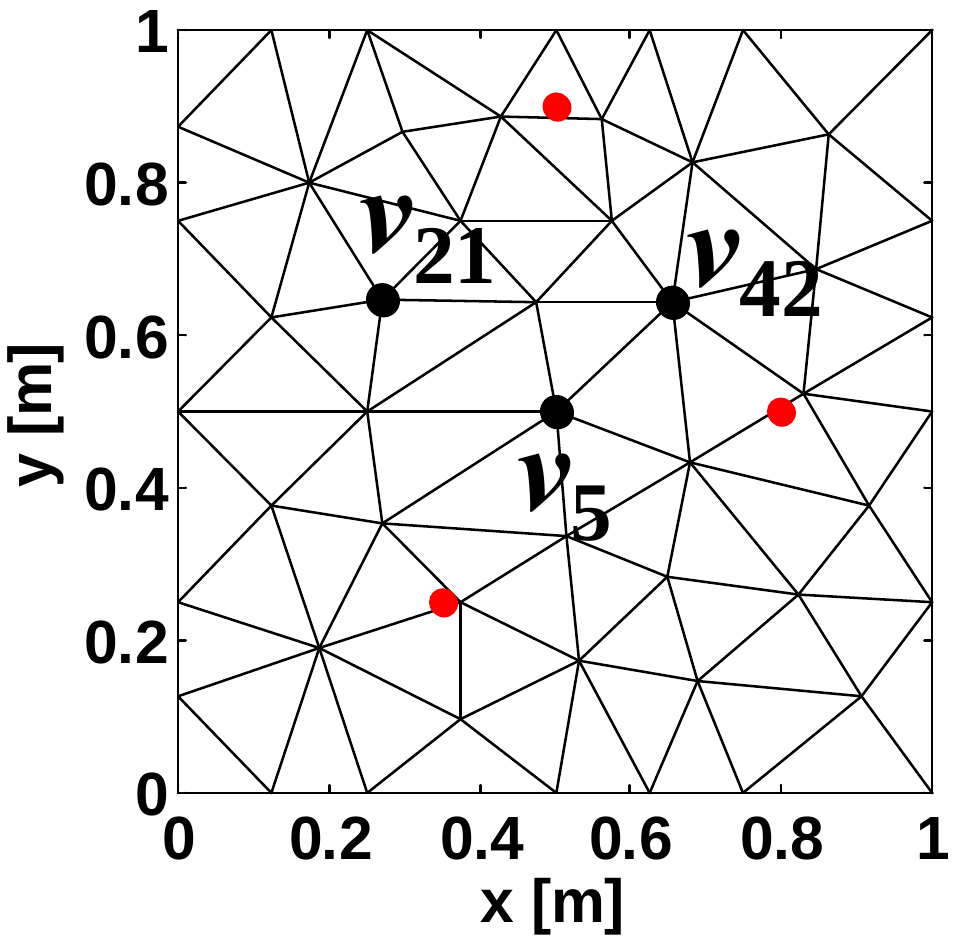}
    }
    \hspace{0.5 cm}
    \subfloat[\label{case2_t50}]{%
      \includegraphics[width=1.8in]{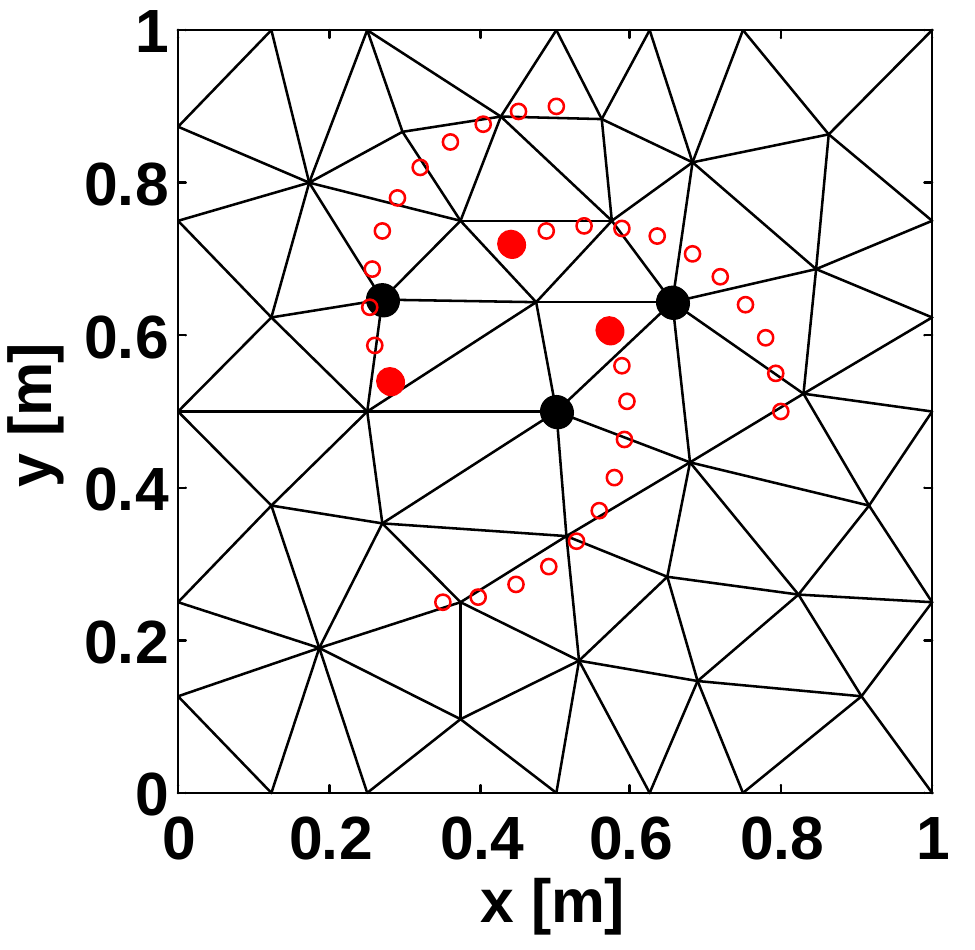}
    }
    \\
   	\subfloat[\label{case2_t100}]{%
      \includegraphics[width=1.8in]{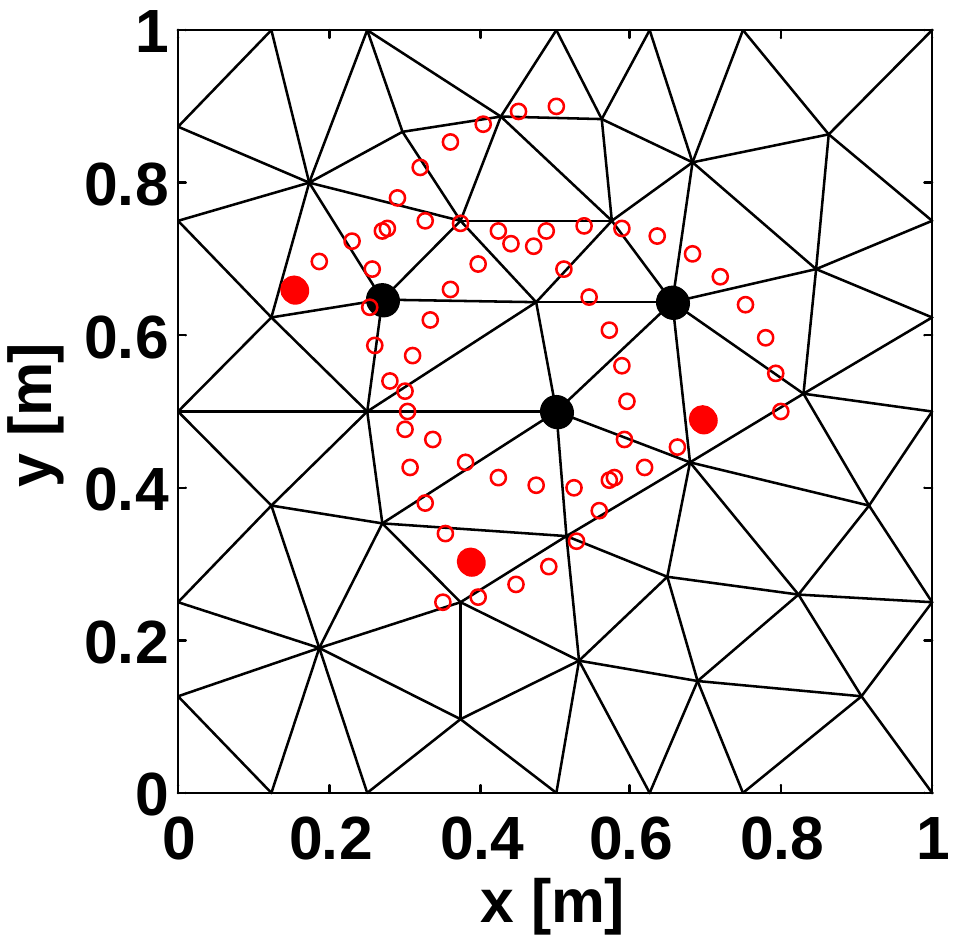}
    }
    \hspace{0.5 cm}
    \subfloat[\label{case2_t200}]{%
      \includegraphics[width=1.8in]{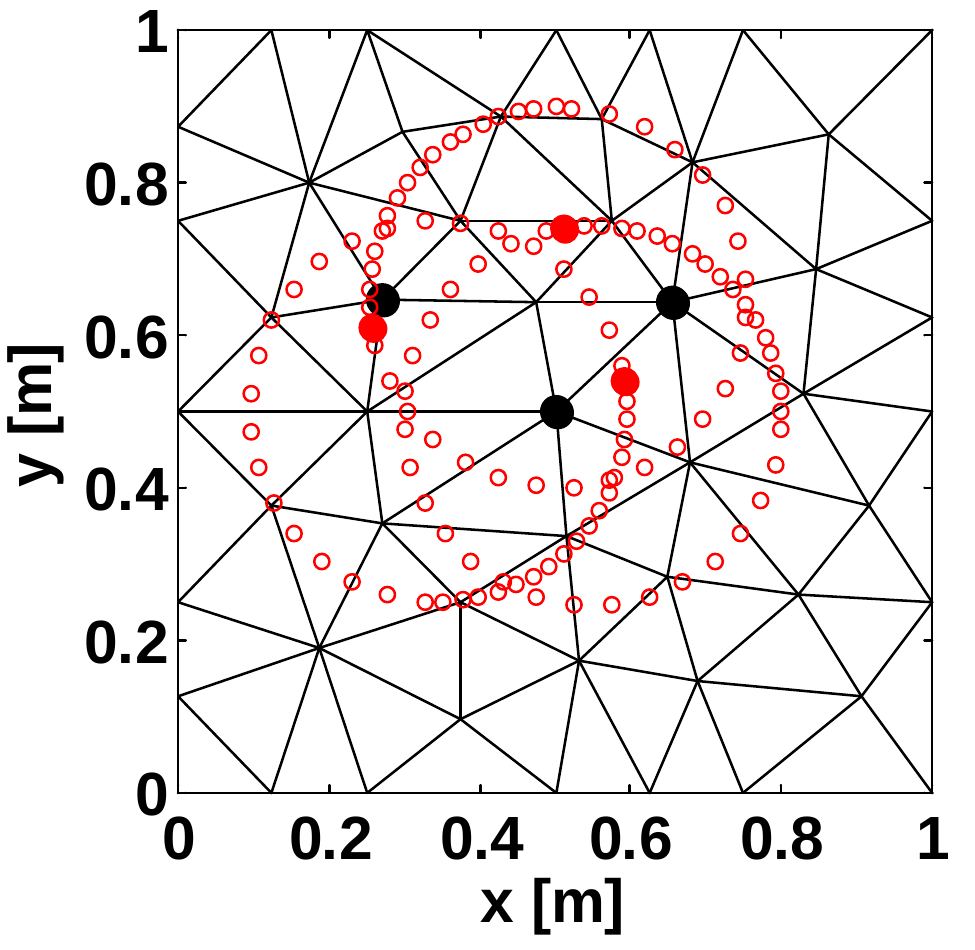}
    }
    \caption{Movement of three particles in the uniform static magnetic field at different time instants ($\Delta t$ = 0.1 ns): (a) t = 0, (b) t = 50$\Delta t$, (c) t = 100$\Delta t$, and (d) t = 200$\Delta t$.}
    \label{particle.move2}
\end{figure}

As second example, Fig. \ref{particle.move2} shows the movement of three negatively-charged particles at different time steps. Similarly as before, these particles describe circular motions because the influence of the static magnetic field is more dominant than interactions among the particles. Particles with positive charges, which are not shown in Fig. \ref{particle.move2}, are again stationary at all time steps due to their much larger masses. We select three random vertices $\nu_{5}$, $\nu_{21}$, $\nu_{42}$ as illustrated in Fig. \ref{case2_t0} for the verification of Gauss' law.
% begin red-colored text
The discrete version of Gauss' law at $t=n\Delta t$, i.e.,
$\widetilde{\mathbf{S}}\cdot \left[\star_{\epsilon}\right]\cdot \mathbf{e}^n = \mathbf{q}^n$
is computed in double-precision floating-point arithmetic. Table \ref{T.Gauss.verify} shows the left- and right-hand side values of this equation and the residuals at several time steps up to $10^6$. The agreement is excellent, and includes at least thirteenth significant digits in all cases and relatively negligible residuals at very large time steps. Note that Gauss' law for magnetism $\mathbf{S}\cdot\mathbf{b}^n=0$ is trivially preserved because only the $B_z$ component is present in this case, and is invariant with respect to $z$.
% end red-colored text

\begin{table}[t]
\begin{center}
\renewcommand{\arraystretch}{1.3}
\caption{Verification of the discrete Gauss' law at different time steps and (global) vertices.}
\begin{tabular}{ccccc}
    \hline
    Vertex & $n$ & $\widetilde{\mathbf{S}}\cdot \left[\star_{\epsilon}\right]\cdot \mathbf{e}^n$ & $\mathbf{q}^n$ & $\widetilde{\mathbf{S}}\cdot \left[\star_{\epsilon}\right]\cdot \mathbf{e}^n - \mathbf{q}^n$\\
    \hline
    \multirow{6}{*}{$\nu_5$}    & $10^1$ &
									-6.206610172341678 $\times 10^{-36}$ &
									                                   0 &
									-6.206610172341678 $\times 10^{-36}$ \\
                                & $10^2$ &
									-3.655787431057314 $\times 10^{-34}$ &
																	   0 &
									-3.655787431057314 $\times 10^{-34}$ \\
                                & $10^3$ &
									-3.996030839009677 $\times 10^{-20}$ &
									-3.996030839009684 $\times 10^{-20}$ &
									 6.620384183831123 $\times 10^{-35}$ \\
                                & $10^4$ &
									-3.581126715387582 $\times 10^{-20}$ &
									-3.581126715385507 $\times 10^{-20}$ &
									-2.074587661969626 $\times 10^{-32}$ \\
                                & $10^5$ &
									 1.442950348685220 $\times 10^{-31}$ &
																	   0 &
									 1.442950348685220 $\times 10^{-31}$ \\
								& $10^6$ &
									-2.830713776667131 $\times 10^{-30}$ &
																	   0 &
									-2.830713776667131 $\times 10^{-30}$ \\
    \hline
    \multirow{6}{*}{$\nu_{21}$} & $10^1$ &
									 3.385423730368188 $\times 10^{-36}$ &
																	   0 &
									 3.385423730368188 $\times 10^{-36}$ \\
                                & $10^2$ &
									-2.045237819570873 $\times 10^{-20}$ &
									-2.045237819570803 $\times 10^{-20}$ &
									-6.981496048403730 $\times 10^{-34}$ \\
                                & $10^3$ &
									-3.751623558930333 $\times 10^{-20}$ &
									-3.751623558929736 $\times 10^{-20}$ &
									-5.970382827600431 $\times 10^{-33}$ \\
                                & $10^4$ &
									-6.801441860813513 $\times 10^{-20}$ &
									-6.801441860811823 $\times 10^{-20}$ &
									-1.690003526199800 $\times 10^{-32}$ \\
                                & $10^5$ &
									-1.264643842388650 $\times 10^{-31}$ &
																       0 &
									-1.264643842388650 $\times 10^{-31}$ \\
								& $10^6$ &
									-1.345470457792998 $\times 10^{-30}$ &
																	   0 &
									-1.345470457792998 $\times 10^{-30}$ \\									
    \hline
    \multirow{6}{*}{$\nu_{42}$} & $10^1$ &
									-4.988033347936703 $\times 10^{-20}$ &
									-4.988033347936874 $\times 10^{-20}$ &
									 1.709262825643672 $\times 10^{-33}$ \\
                                & $10^2$ &
									-3.159139423392524 $\times 10^{-20}$ &
									-3.159139423392751 $\times 10^{-20}$ &
									 2.268986215731212 $\times 10^{-33}$ \\
                                & $10^3$ &
									-3.915500401400283 $\times 10^{-20}$ &
									-3.915500401401157 $\times 10^{-20}$ &
									 8.738907122657083 $\times 10^{-33}$ \\
                                & $10^4$ &
									 8.443811020826656 $\times 10^{-32}$ &
											                           0 &
									 8.443811020826656 $\times 10^{-32}$ \\
                                & $10^5$ &
									 7.919916698552179 $\times 10^{-31}$ &
													                   0 &
									 7.919916698552179 $\times 10^{-31}$ \\
								& $10^6$ &
									-5.373841942680172 $\times 10^{-21}$ &
									-5.373841950539624 $\times 10^{-21}$ &
									 7.859451526095034 $\times 10^{-30}$ \\		
    \hline
\end{tabular}
\label{T.Gauss.verify}
\end{center}
\end{table}

\begin{figure}[t]
    \centering
    \includegraphics[width=2.5in]{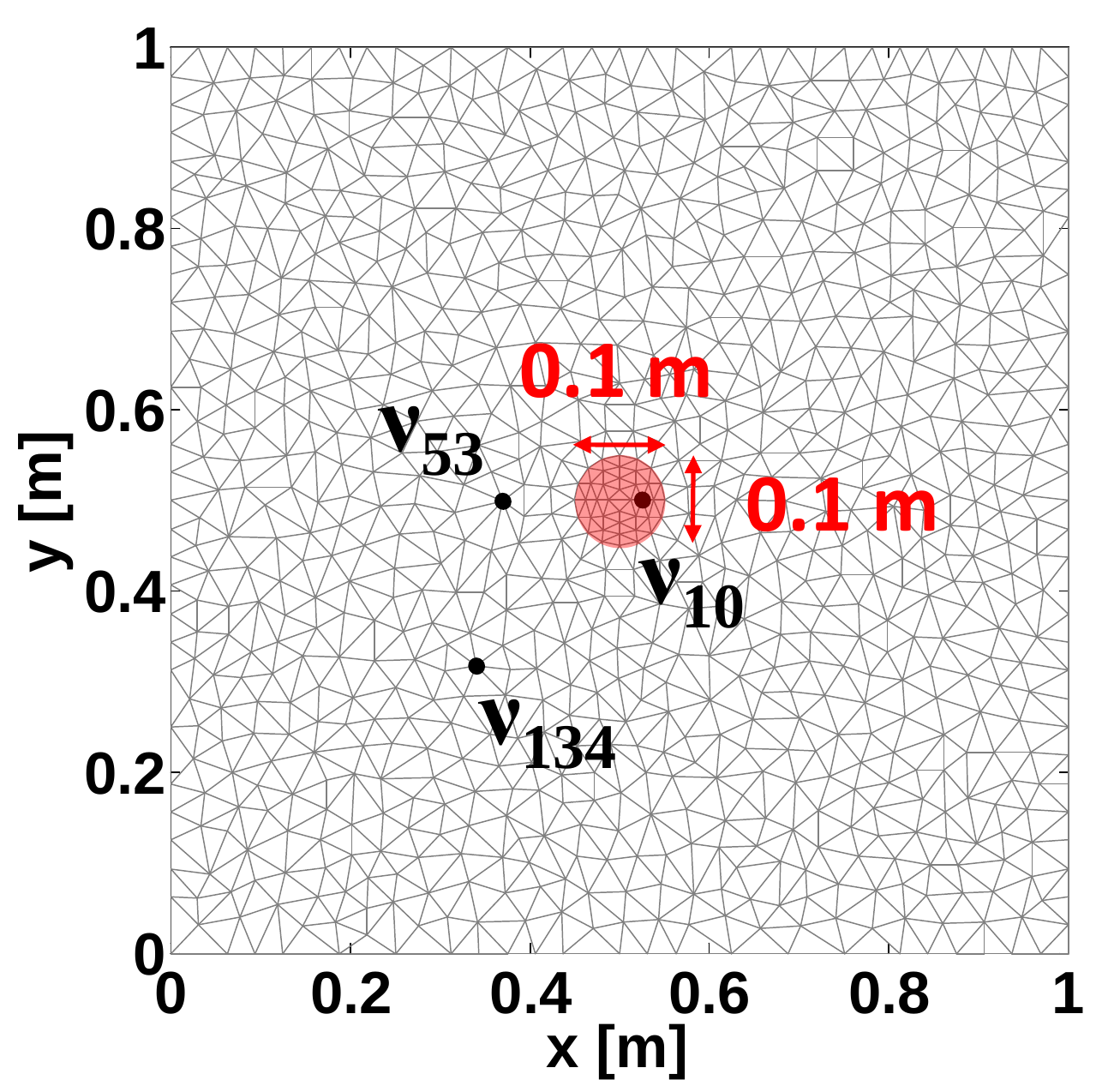}
    \caption{Mesh with 2539 edges and three selected nodes $\nu_{10}$, $\nu_{53}$, and $\nu_{134}$. A total of $4\times 10^3$  negatively charged particles are initially placed in the red circle, uniformily distributed.}
    \label{mesh.DoF2539}
\end{figure}

\begin{figure}[t]
	\centering
	\subfloat[\label{case3_t10k}]{%
      \includegraphics[width=1.8in]{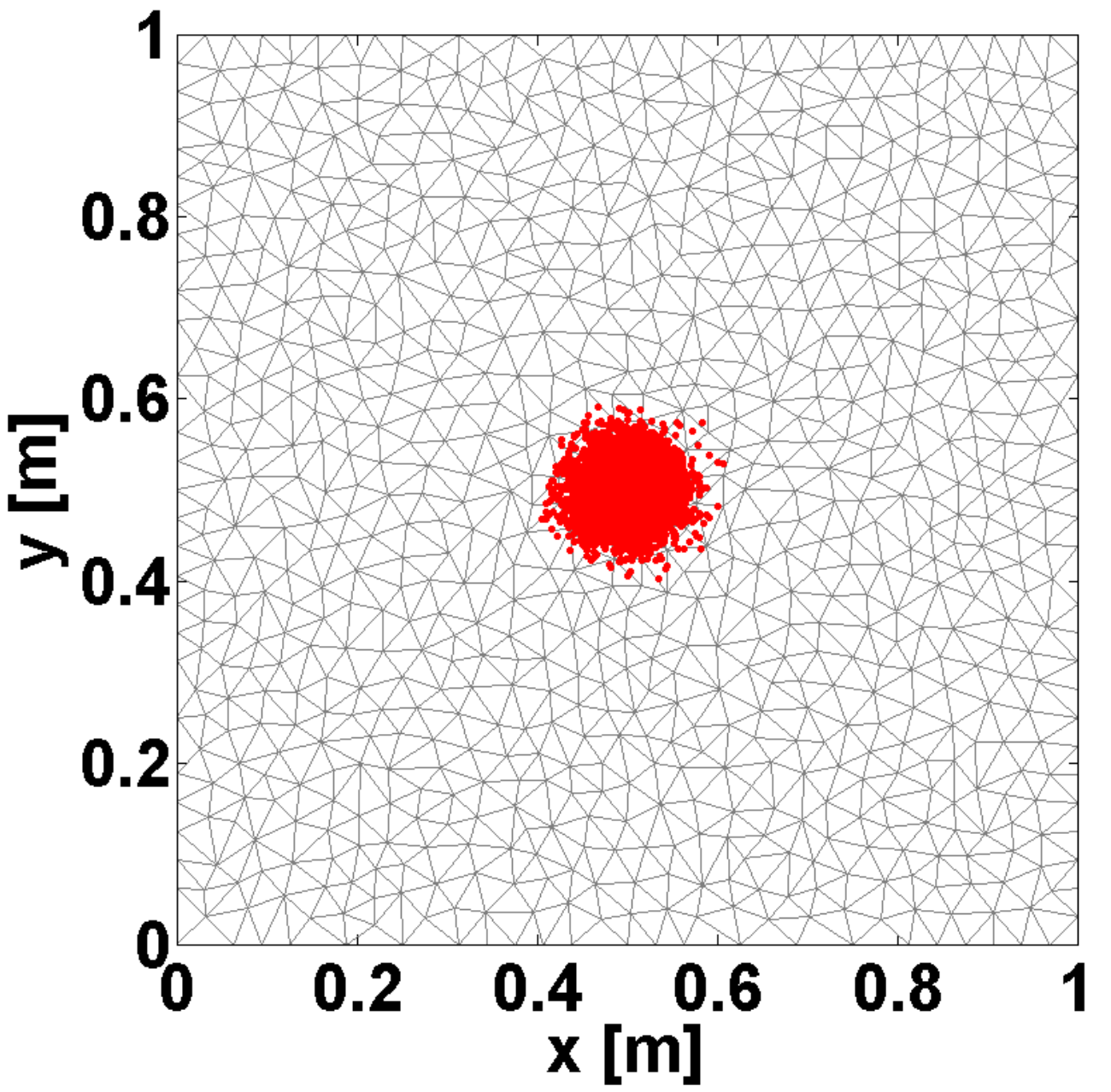}
    }
    \hspace{0.5 cm}
    \subfloat[\label{case3_t20k}]{%
      \includegraphics[width=1.8in]{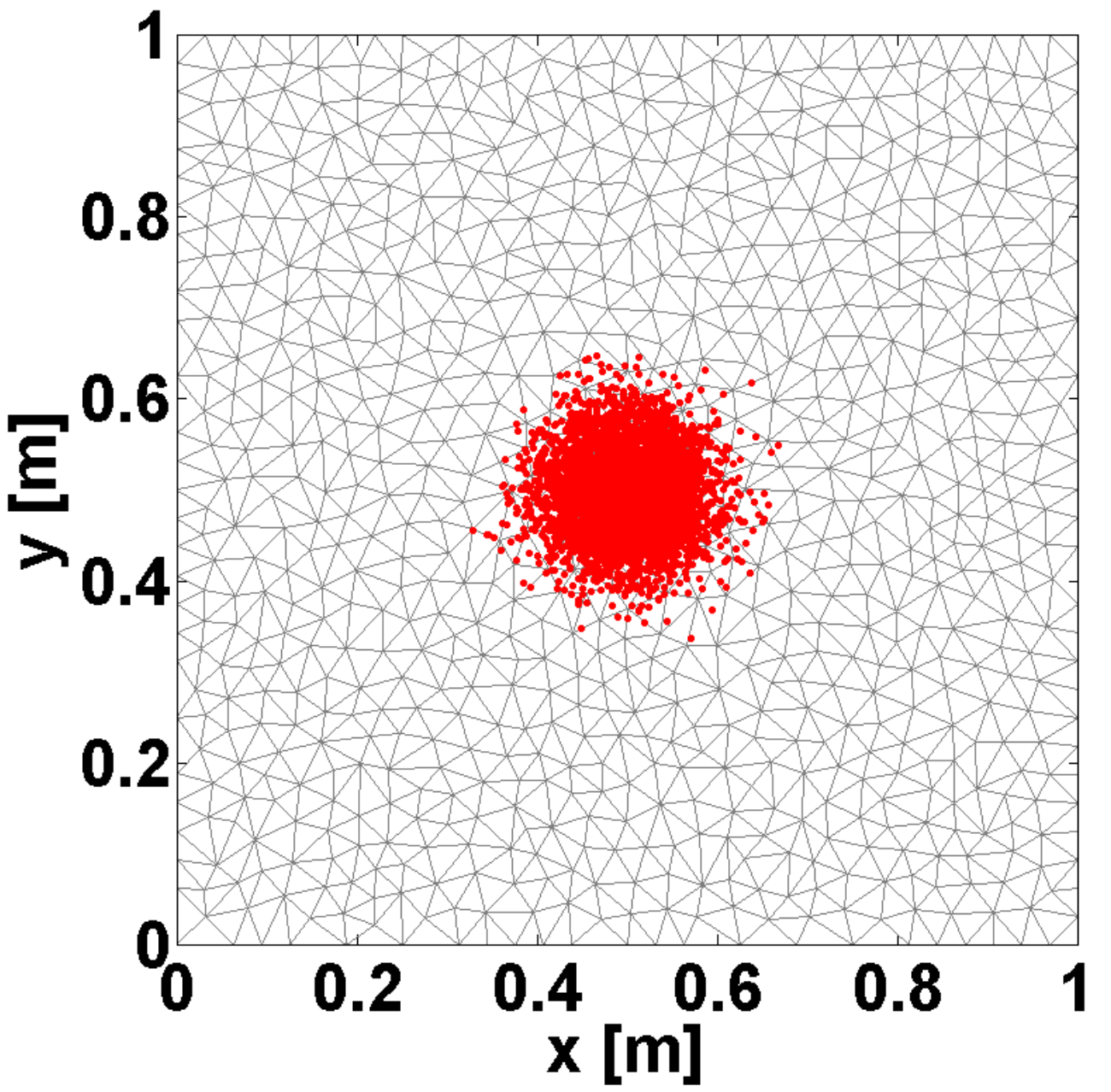}
    }
    \\
   	\subfloat[\label{case3_t40k}]{%
      \includegraphics[width=1.8in]{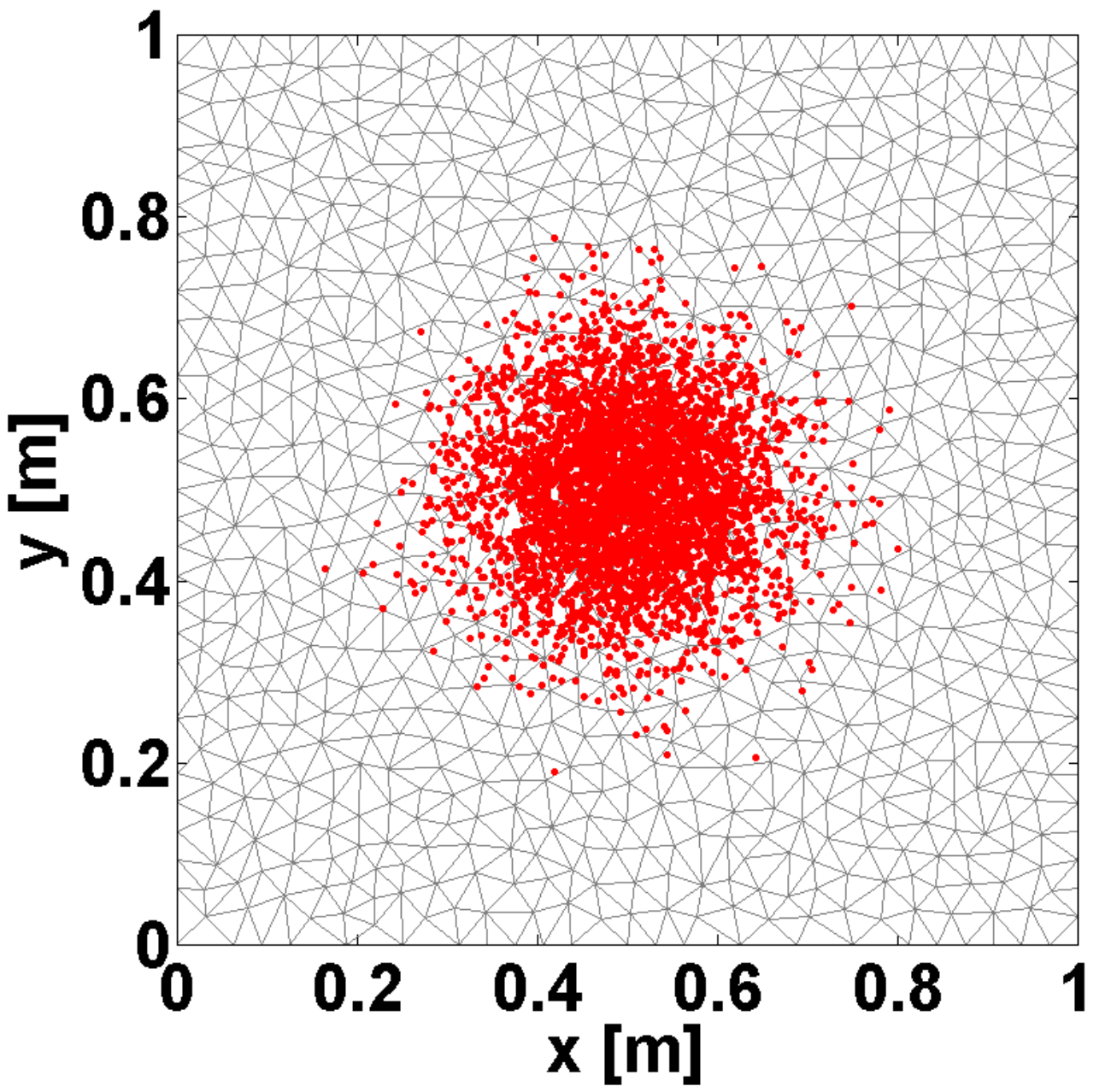}
    }
    \hspace{0.5 cm}
    \subfloat[\label{case3_t60k}]{%
      \includegraphics[width=1.8in]{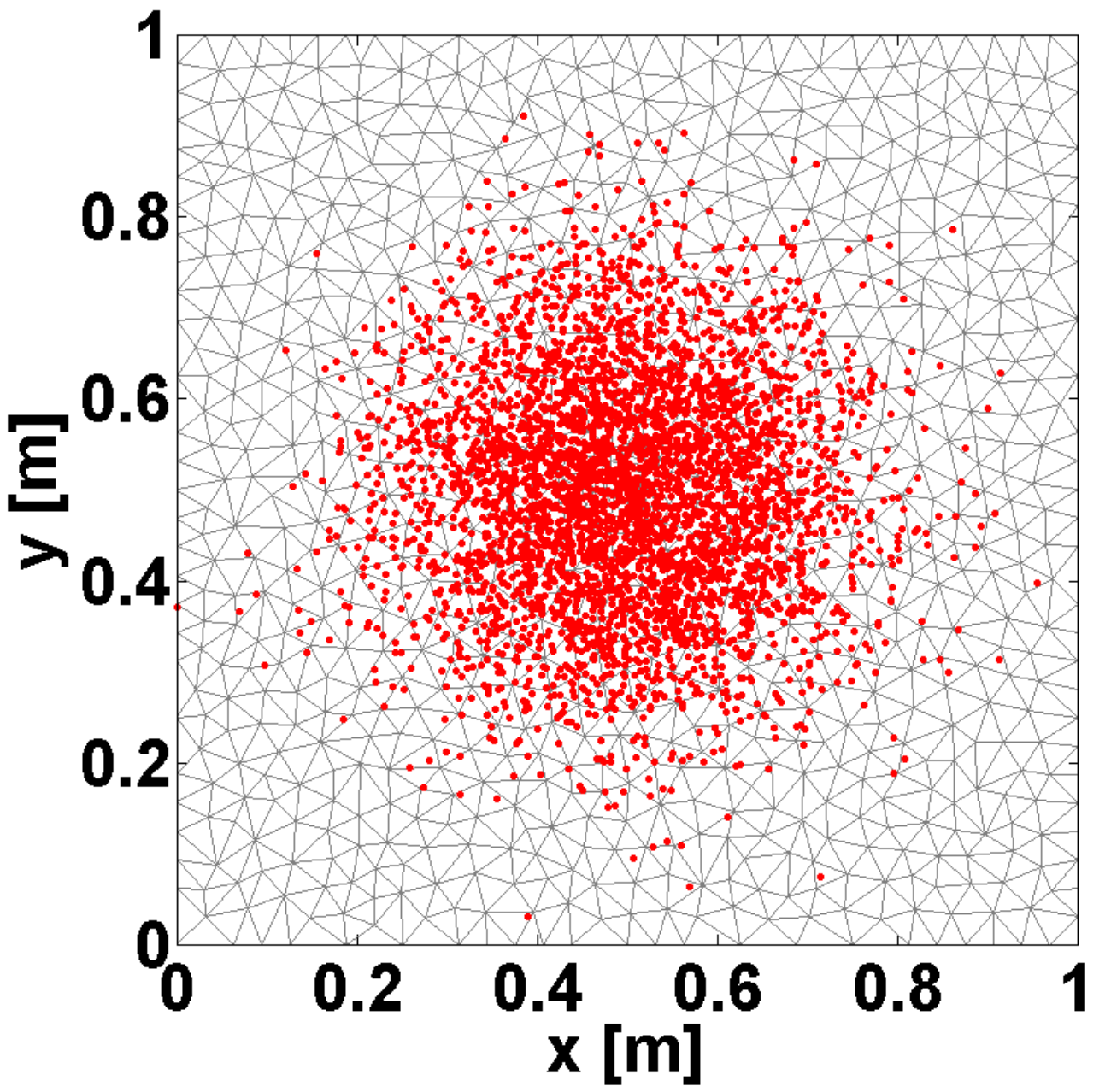}
    }
    \caption{Distribution of $4 \times 10^3$ particles with initial Maxwellian distribution, and zero initial fields. The particle distribution is shown at different time instants ($\Delta t$ = 0.01 ns): (a) $t = 10^4\Delta t$, (b) $t = 2 \times 10^4\Delta t$, (c) $t = 4 \times 10^4\Delta t$, and (d) $t = 6 \times 10^4 \Delta t$.}
    \label{particle.move3}
\end{figure}

\begin{figure}[t]
    \centering
    \includegraphics[width=3.0in]{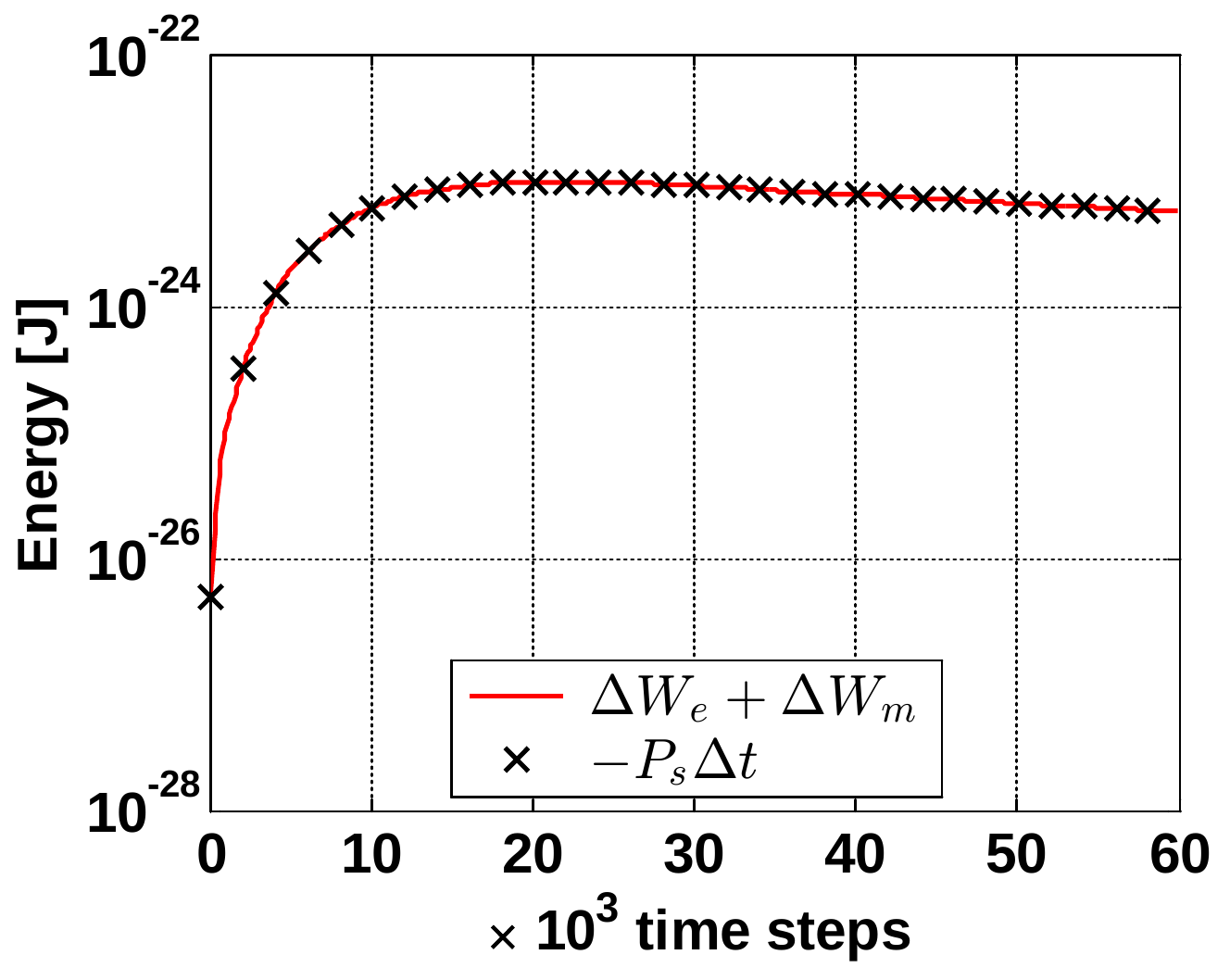}
    \caption{Numerical verification of energy conservation with a plot of the left hand side and the right hand side of \eqref{energy.bal.4} for all time steps.}
    \label{energy.comp}
\end{figure}

% begin blue-colored text
For the third example, we consider the PIC simulation of a blowing-up plasma ball (circle) composed of two species: electrons (hot) and ions (cold). Initially,
all electrons and ions are overlapped, so that the local charge is zero everywhere. The mesh is dense around the initial plasma ball and relaxed radially. The mesh is depicted in Fig. \ref{mesh.DoF2539} and has 2539 edge elements.
 A total of 4000 negatively charged particles are initially randomly distributed inside the red circle shown, with particle density $n_e=4\times 10^3/(0.05^2\pi)=5.0930\times 10^5$ m$^{-3}$. Electron velocities are initialized with a Maxwellian distribution, with a thermal velocity $|\mathbf{v}_{th}| = 10^{-3} c$ m/s, where c is the light speed (nonrelativistic regime). Positive ions are assumed much more massive, with zero velocity.
In Fig. \ref{mesh.DoF2539}, three nodes $\nu_{10}$, $\nu_{53}$, and $\nu_{134}$ are designated for veryfying charge conservation at all times. The electron Debye length is such that $\lam_D^2=\epsilon_0 k T/(n_e)^{3/2}q^2$, which gives $\lam_D=0.1974$ m from the settings above. 
Fig. \ref{particle.move3} shows the distribution of the 4000 particles at different time steps, illustrating the expansion of the plasma ball.
% end blue-colored text
% begin blue-colored text
To examine energy conservation, we consider the energy balanced equation
\begin{flalign}
\frac{\pa}{\pa t}
	\left(\frac{1}{2}\mathbf{E}\cdot\epsilon\mathbf{E} + \frac{1}{2}\mathbf{B}\cdot\mu^{-1}\mathbf{B}
	\right) + \mathbf{E}\cdot\mathbf{J} = 0. \label{energy.bal}
\end{flalign}
After spatial discretization, \eqref{energy.bal} writes 
\begin{flalign}
\frac{d}{dt}
	\left(\frac{1}{2}\mathbf{e}^T\cdot\left[\star_{\epsilon}\right]\cdot\mathbf{e}
		+ \frac{1}{2}\mathbf{b}^T\cdot\left[\star_{\mu^{-1}}\right]\cdot\mathbf{b}
	\right) + \mathbf{e}^T\cdot\mathbf{i} = 0, \label{energy.bal.2}
\end{flalign}
or, more concisely,
\begin{flalign}
\frac{d}{dt}
	\left(W_e + W_m\right) + P_s = 0, \label{energy.bal.3}
\end{flalign}
where $W_e$ and $W_m$ are the electric and magnetic energy density terms, and $P_s$ is the term associated with the presence of electric current $\mathbf{J}$ from the moving charges. Using a leap-frog scheme for time-discretization, we obtain
\begin{flalign}
\Delta W_e^{n+\frac{1}{2}} + \Delta W_m^{n+\frac{1}{2}} =  -P_s^{n+\frac{1}{2}}\Delta t, \label{energy.bal.4}
\end{flalign}
where half-integer times are considered to coincide with $\mathbf{i}$. Figure \ref{energy.comp} shows the comparison between the left hand side and the right hand side of \eqref{energy.bal.4} for all time steps. An excellent agreement is observed, which numerically verifies energy conservation.
% end blue-colored text

% begin blue-colored text
The discrete version of Gauss' law is also examined for this case in Table \ref{T.Gauss.verify.2}. Again, there is a very good match between the two terms of Gauss' law for all times, with at least eight significant digits of agreement even after a million time steps.
% end blue-colored text
\begin{table}[t]
\begin{center}
\renewcommand{\arraystretch}{1.3}
\caption{Verification of the discrete Gauss' law for PIC simulations with many particles at different time steps and for three arbitrary (global) vertices.}
\begin{tabular}{ccccc}
    \hline
    Vertex & $n$ & $\widetilde{\mathbf{S}}\cdot \left[\star_{\epsilon}\right]\cdot \mathbf{e}^n$ & $\mathbf{q}^n$ & $\widetilde{\mathbf{S}}\cdot \left[\star_{\epsilon}\right]\cdot \mathbf{e}^n - \mathbf{q}^n$\\
    \hline
    \multirow{5}{*}{$\nu_{10}$} & $10^1$ &
									 3.938626419217293 $\times 10^{-21}$ &
									 3.938828609135238 $\times 10^{-21}$ &
									-2.021899179450384 $\times 10^{-25}$ \\
                                & $10^2$ &
									 4.216965415302240 $\times 10^{-20}$ &
									 4.216619000000000 $\times 10^{-20}$ &
									 3.464153022399240 $\times 10^{-24}$ \\
                                & $10^3$ &
									 5.693133763020551 $\times 10^{-19}$ &
									 5.692525000000000 $\times 10^{-19}$ &
									 6.087630205511676 $\times 10^{-23}$ \\
                                & $10^4$ &
									 3.376346144358901 $\times 10^{-18}$ &
									 3.373843000000000 $\times 10^{-18}$ &
									 2.503144358901224 $\times 10^{-21}$ \\
                                & $6\times10^4$ &
									 1.948941527506024 $\times 10^{-17}$ &
									 1.947695000000000 $\times 10^{-17}$ &
									 1.246527506023833 $\times 10^{-20}$ \\
    \hline
    \multirow{6}{*}{$\nu_{53}$} & $10^1$ &
									-6.035730222106830 $\times 10^{-25}$ &
																	   0 &
									-6.035730222106830 $\times 10^{-25}$ \\
                                & $10^2$ &
									-3.859548532081658 $\times 10^{-24}$ &
																	   0 &
									-3.859548532081658 $\times 10^{-24}$ \\
                                & $10^3$ &
									-3.778195570243296 $\times 10^{-23}$ &
																	   0 &
									-3.778195570243296 $\times 10^{-23}$ \\
                                & $10^4$ &
									-2.171627281773591 $\times 10^{-21}$ &
																	   0 &
									-2.171627281773591 $\times 10^{-21}$ \\
                                & $6\times10^4$ &
									-6.202694078733229 $\times 10^{-18}$ &
									-6.183721000000000 $\times 10^{-18}$ &
									-1.897307873322922 $\times 10^{-20}$ \\
    \hline
    \multirow{6}{*}{$\nu_{134}$}& $10^1$ &
									4.216746669112738 $\times 10^{-31}$ &
																	  0 &
									4.216746669112738 $\times 10^{-31}$ \\
                                & $10^2$ &
									 1.762689521069704 $\times 10^{-26}$ &
																	   0 &
									 1.762689521069704 $\times 10^{-26}$ \\
                                & $10^3$ &
									 6.830106418522705 $\times 10^{-26}$ &
																	   0 &
									 6.830106418522705 $\times 10^{-26}$ \\
                                & $10^4$ &
									 1.532866351101663 $\times 10^{-24}$ &
											                           0 &
									 1.532866351101663 $\times 10^{-24}$ \\
                                & $6\times10^4$ &
									-1.478546530191113 $\times 10^{-18}$ &
									-1.478877000000000 $\times 10^{-18}$ &
									 3.304698088865041 $\times 10^{-22}$ \\							
    \hline
\end{tabular}
\label{T.Gauss.verify.2}
\end{center}
\end{table}

%%%%%%%%%%%%%%%%%%%%%%%%%%%%%%%%%%%%%%%%%%%%%%%%%%%%%%%%%%%%%%%%%%%%%%%%%%%%%%%%%%%%%%%%%%%%%%%%%%%%%%%%%%%%%%
\section{Concluding Remarks}
%%%%%%%%%%%%%%%%%%%%%%%%%%%%%%%%%%%%%%%%%%%%%%%%%%%%%%%%%%%%%%%%%%%%%%%%%%%%%%%%%%%%%%%%%%%%%%%%%%%%%%%%%%%%%%
A new, geometrically intuitive charge-conserving scatter-gather algorithm for full electromagnetic PIC simulations has been presented for arbitrary unstructured grids. The algorithm relies upon the representation of the various dynamical quantities as discrete differential forms of different degrees, and on their self-consistent interpolation by Whitney forms. Preservation of Gauss' law is demonstrated for all times, both analytically and by means of numerical tests.

\section*{Acknowledgments}
This work was supported in part by NSF under grant ECCS-1305838 and OSC under grants PAS-0061 and PAS-0110. The authors would like to thank the reviewers for pointing out reference~\cite{Pinto14:Charge} and for their suggested clarifications to the text.

\appendix
%%%%%%%%%%%%%%%%%%%%%%%%%%%%%%%%%%%%%%%%%%%%%%%%%%%%%%%%%%%%%%%%%%%%%%%%%%%%%%%%%%%%%%%%%%%%%%%%%%%%%%%%%%%%%%
\section{Whitney forms: Basic properties}
\label{app.a}
%%%%%%%%%%%%%%%%%%%%%%%%%%%%%%%%%%%%%%%%%%%%%%%%%%%%%%%%%%%%%%%%%%%%%%%%%%%%%%%%%%%%%%%%%%%%%%%%%%%%%%%%%%%%%%
For convenience, we provide here the explicit expressions of Whitney forms~\cite{Whitney:Geometric} in 3-D. 
In the past, Whitney forms have proved useful in finite element modeling of electromagnetic fields~\cite{Bossavit:Computational, Jin:Finite, Bondeson:Computational}, to suppress spurious modes. Although Whitney forms can be more succinctly and elegantly expressed using the exterior calculus of differential forms~\cite{Kim11:Parallel, Teixeira99:Lattice, He07:Differential, Whitney:Geometric}, we adopt here the more familiar notation of vector calculus.

\begin{figure}[t]
	\centering
	\subfloat[\label{bary.line}]{%
      \includegraphics[width=2.0in]{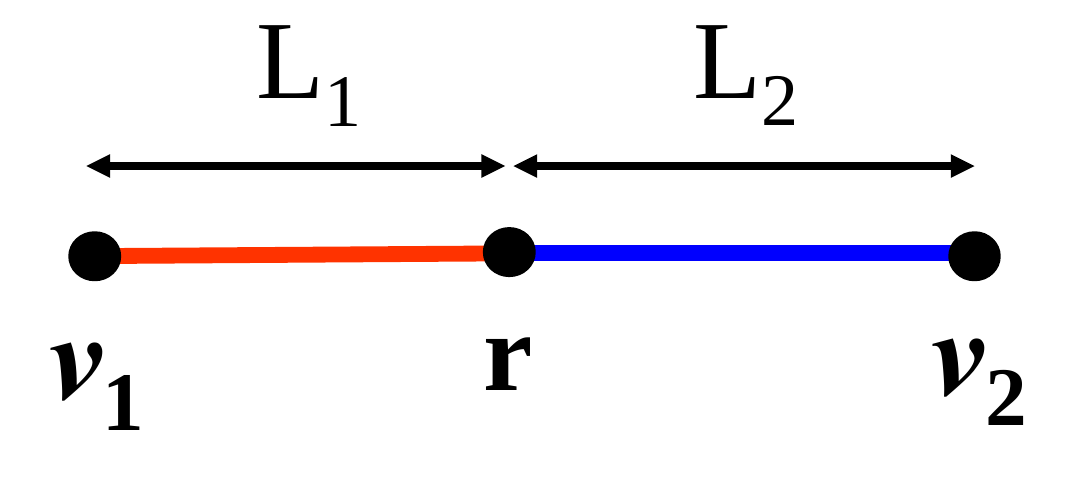}
    }
    \hspace{0.5 cm}
    \subfloat[\label{bary.tri}]{%
      \includegraphics[width=2.4in]{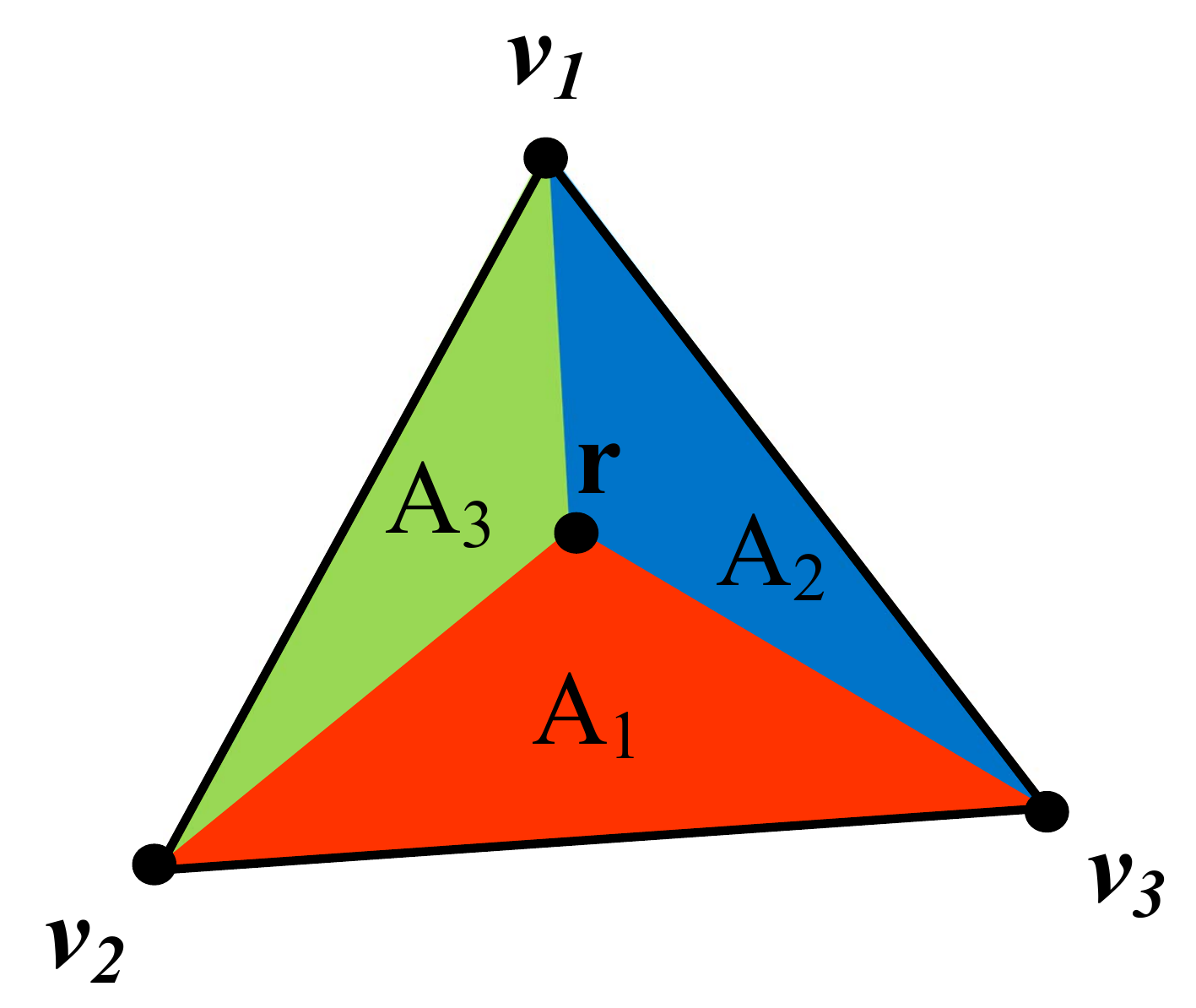}
    }
    \caption{Geometric illustration for Whitney 0-forms (barycentric coordinates) of a point $\mathbf{r}$ in simplices of various degrees: (a) 1-D simplex and (b) 2-D simplex.}
    \label{bary.coord}
\end{figure}

In 3-D, there are four types of Whitney $p$-forms, according to their degree $p$.
A {\it Whitney 0-form} is a continuous scalar function simply expressed as~\cite{Teixeira99:Lattice}
\begin{flalign}
W_i^0 (\mathbf{r})= \lam_i(\mathbf{r}), \label{Whitney.0form}
\end{flalign}
where the subscript $i$ represents vertex $i$ and $\lam_i$ is the barycentric coordinate~\cite{Silvester:Finite} associated with vertex $i$. The geometric construction for barycentric coordinates is illustrated in Fig. \ref{bary.coord}. For a 1-D simplex (i.e. edge), the barycentric coordinates associated to the vertices $\nu_1$ and $\nu_2$ of any point $\mathbf{r}$ in the simplex  are equal to ratios $\lam_1=L_1/(L_1+L_2)$ and $\lam_2=L_2/(L_1+L_2)$, respectively, with $L_1$ and $L_2$ as indicated in Fig. \ref{bary.coord}. For a 2-D simplex (triangle), the barycentric coordinates associated to the three vertices $\nu_1$, $\nu_2$ and $\nu_3$ of any point $\mathbf{r}$ in the simple are equal to $\lam_1=A_1/A$, $\lam_2=A_2/A$, and $\lam_3=A_3/A$, respectively, with the areas $A_1$, $A_2$, and $A_3$ as indicated and $A=A_1+A_2+A_3$. In a 3-D simplex, which is a tetrahedron, the barycentric coordinates can be similarly written as volume ratios. It is clear that $0 \leq \lam_i \leq 1$ for all $i$ and that the sum of the barycentric coordinates of any given point $\mathbf{r}$ associated to the neighbor vertices equals to one. Hereinafter, the dependence on $\mathbf{r}$ is dropped for notational simplicity, i.e., $\lam_i(\mathbf{r})=\lam_i$.

\begin{figure}[t]
	\centering
	\subfloat[\label{Whitney_1form_L}]{%
      \includegraphics[width=2.4in]{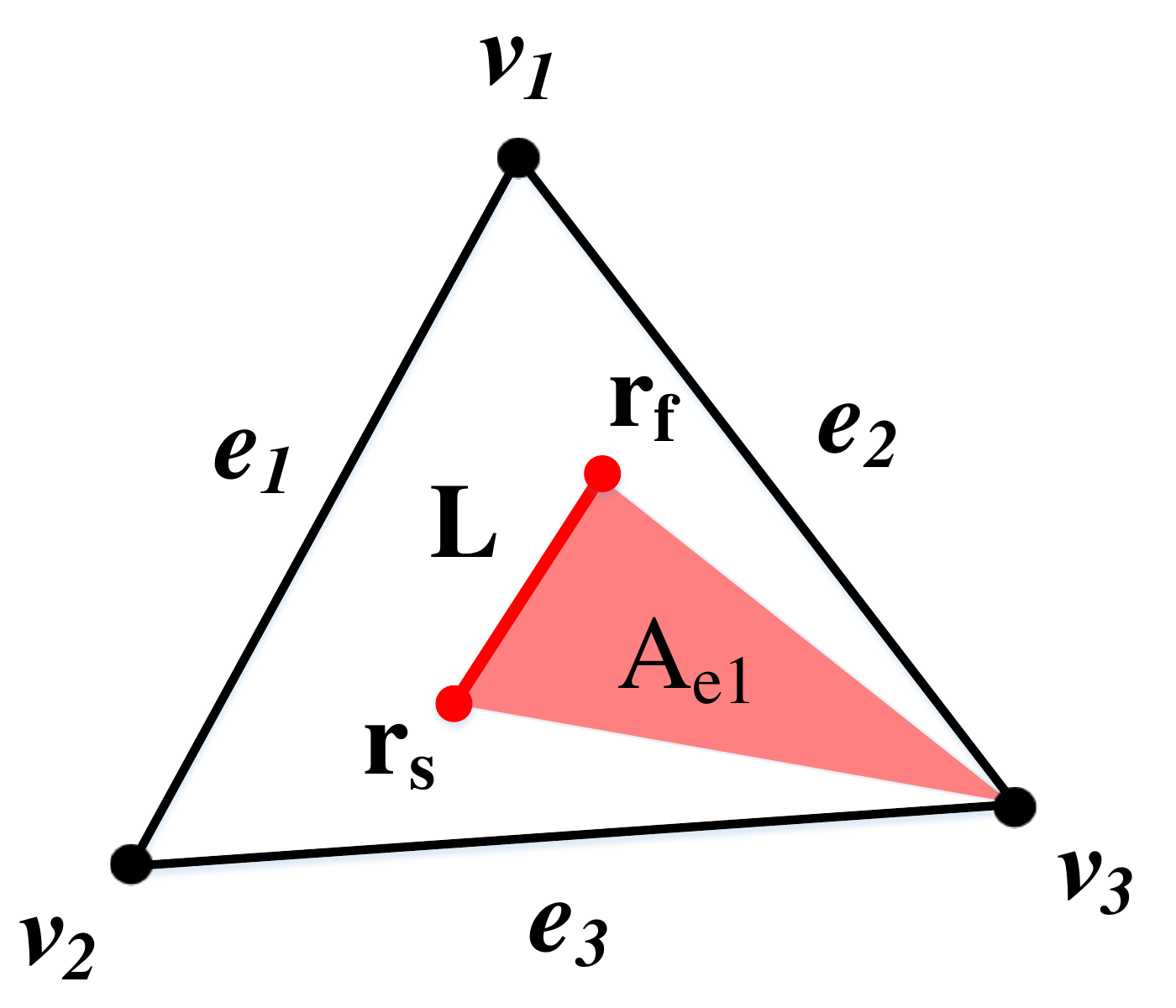}
    }
    \hspace{0.5 cm}
    \subfloat[\label{Whitney_1form_dL}]{%
      \includegraphics[width=2.4in]{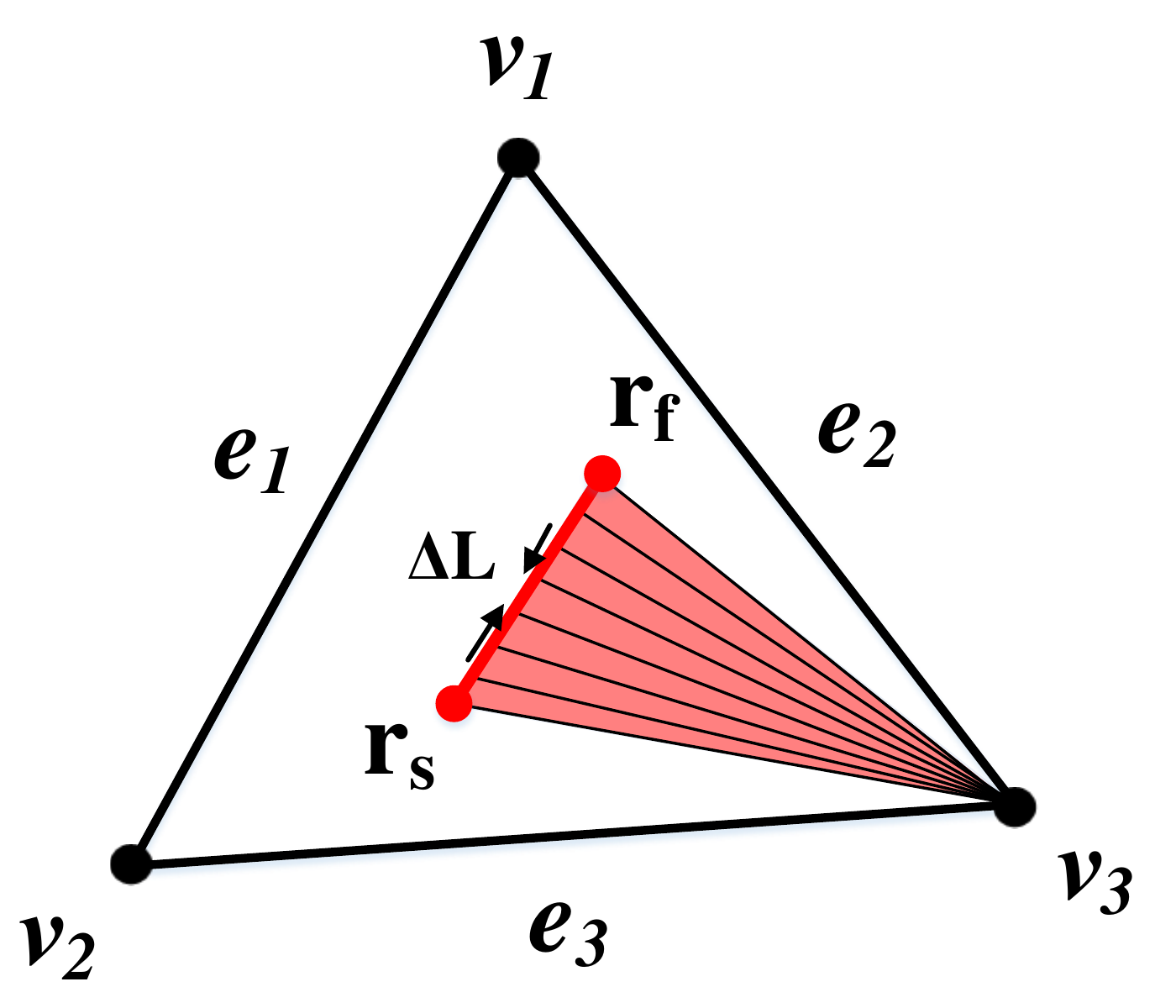}
    }
    \caption{Geometric illustration of the weight assigned to Whitney 1-forms representing a segment $\mathbf{L}$ in a 2-D simplex: (a) In red color is the area $A_{e1}$ associated with the Whitney 1-form on $e_1$ (edge 1) that represents $\mathbf{L}$. The associated weight is given by $A_{e1}/A$,
    where $A$ is the total area of the triangle composed of $\nu_1, \nu_2$, and $\nu_3$. A similar construction can be made for the other two edges $e_2$ and $e_3$. (b) Area 
    represented by a sum of small triangles. See the main text for more details.}
    \label{Whitney.1form.area}
\end{figure}

The vector (function) proxy of a {\it Whitney 1-form} associated with an arbitrary edge $ij$\footnote{For the sake of clarity, we adopt in this Appendix a vertex-based indexing for all types of elements. This is in contrast to the single-indexing adopted for all element types elsewhere in the paper.} bounded by vertices $i$ and $j$ is expressed as~\cite{Teixeira99:Lattice}
\begin{flalign}
\mathbf{W}_{ij}^1 (\mathbf{r})= \lam_i\na\lam_j - \lam_j \na\lam_i. \label{Whitney.1form}
\end{flalign}
% begin red-colored text
For a brief geometric illustration of the Whitney 1-form, let us consider Fig.~\ref{Whitney.1form.area}. The area $A_{e1}$, which is associated with $e_1$, is
\begin{flalign}
A_{e1}  = A \left[\lam_1^s\lam_2^f - \lam_2^s\lam_1^f \right],
\end{flalign}
where $\lambda_i^s$ and $\lambda_i^f$ are shorthands of $\lambda_i(\mathbf{r}_s)$ and $\lambda_i(\mathbf{r}_f)$, respectively. As Fig.~\ref{Whitney_1form_dL} shows, $A_{e1}$ can be regarded as the sum of small triangles such that
\begin{flalign}
A_{e1} & = A \sum_n\left[\lam_1^n\lam_2^{n+1} - \lam_2^n\lam_1^{n+1} \right] \notag\\
	& = A \sum_n\left[\lam_1^n(\lam_2^n+\Delta\lam_2^n) - \lam_2^n(\lam_1^n+\Delta\lam_1^n) \right] \notag\\
	& = A \sum_n\left[\lam_1^n \Delta\lam_2^n - \lam_2^n \Delta\lam_1^n \right]. \label{Whitney.1form.area.a}
\end{flalign}
After taking the limit of infinitesimally small triangles and transforming this summation to an integral, we obtain
\begin{flalign}
A_{e1}
= A \int_{\mathbf{r}_s}^{\mathbf{r}_f}\left[\lam_1 \na\lam_2 - \lam_2 \na\lam_1 \right]\cdot d{\mathbf{L}}
= A \int_{\mathbf{r}_s}^{\mathbf{r}_f} \mathbf{W}_{12}^1 (\mathbf{r}) \cdot d{\mathbf{L}}.
\label{Whitney.1form.area.b}
\end{flalign}
The areas associated with $e_2$ and $e_3$ can be derived in a similar fashion.
The last integral above can be viewed as the generalization of the concept of barycentric coordinates from 0-dimensional objects (points) to 
1-dimensional objects (segments). That is, this relation illustrates that, in the same manner as the Whitney 0-forms (barycentric coordinates) are used to represent a point as a weighted sum 
of nearby vertices $i=1,2,3$ (with respective weights $A_i/A$),
Whitney 1-forms represent {\it any segment} $[\mathbf{r}_s,\mathbf{r}_f]$ {\it in terms of the nearby edges} $e_1$, $e_2$, and $e_3$ (now with weights $A_{e1}/A$, $A_{e2}/A$, and $A_{e3}/A$, respectively).
In both cases, the weights are computed by the ``contraction''~\cite{Teixeira99:Lattice} of the Whitney form with the corresponding geometric object. For a 0-form, this contraction simply means an evaluation of 
$W_i^0$
at the point $\mathbf{r}$, i.e.,
$W_i^0(\mathbf{r})$ as in \eqref{Whitney.0form}, whereas for a 1-form, this contraction means an evaluation of the line integral of $\mathbf{W}_{ij}^1$ along
the segment $[\mathbf{r}_s,\mathbf{r}_f]$ as in \eqref{Whitney.1form.area.b}.  For a more general description of these Whitney form properties, see~\cite{Teixeira14:Lattice}. A comprehensive discussion of the integral of Whitney 1-forms along a straight segment is presented in \ref{app.b} below.

% end red-colored text

Likewise, the vector proxy of a {\it Whitney 2-form} associated with a triangular cell $ijk$ is a vector function expressed as~\cite{Teixeira99:Lattice}
\begin{flalign}
\mathbf{W}_{ijk}^2 (\mathbf{r})= 2
	\Big[
		\lam_i\na\lam_j\times\na\lam_k + \lam_j\na\lam_k\times\na\lam_i + \lam_k\na\lam_i\times\na\lam_j
	\Big]. \label{Whitney.2form}
\end{flalign}
Finally, in 3-D, the proxy of a {\it Whitney 3-form} associated with a tetrahedral cell $ijkl$ in 3-D is a scalar function written as~\cite{Teixeira99:Lattice}
\begin{flalign}
W_{ijkl}^3 (\mathbf{r})&= 6
	\Big[
		\lam_i\na\lam_j\cdot\left(\na\lam_k\times\na\lam_l\right) + 	
		\lam_j\na\lam_k\cdot\left(\na\lam_l\times\na\lam_i\right)	\notag\\
	&\qquad +
		\lam_k\na\lam_l\cdot\left(\na\lam_i\times\na\lam_j\right) +
		\lam_l\na\lam_i\cdot\left(\na\lam_j\times\na\lam_k\right)
	\Big], \label{Whitney.3form}
\end{flalign}
Despite the complicated-looking expression~\ref{Whitney.3form}, $W_{ijkl}^3$ can be shown in 3-D to be simply equal to
\begin{flalign}
W_{ijkl}^3 (\mathbf{r})=
	\begin{cases}
	\frac{1}{V}, \text{ if } \mathbf{r} \text{ is in the tetrahedron $ijkl$},
	\\ 0, \text{ otherwise},
	\end{cases}
\end{flalign}
where $V$ is the volume of the tetrahedron $ijkl$~\cite{Teixeira99:Lattice}.
Whitney forms are interpolatory in the precise sense that they are equal to one when ``evaluated on'' the respective elements (vertices, edges, triangles, and tetrahedra) and to zero on all remaining elements of the grid, where ``evaluated on'' in the case of $\mathbf{W}_{ij}^1$, $\mathbf{W}_{ijk}^2$, and ${W}_{ijkl}^3$ means ``integrated over'' edges, triangles, or tetrahedrons respectively\footnote{That is, line, surface, or volume integration, for a Whitney form of degree $p=$ 1, 2, and 3, respectively.}.
Furthermore, Whitney forms inherit the same type of continuity of the fields they represent.
Specifically, $W_i^0 (\mathbf{r})$ is a continuous scalar function (representing scalar potentials for example), $\mathbf{W}_{ij}^1 (\mathbf{r})$ is a tangentially continuous vector function (representing ``intensity'' vector fields for example), $\mathbf{W}_{ijk}^2 (\mathbf{r})$ is a normally continuous vector functions (representing ``flux density'' vector fields or volumetric current densities for example) and $W_{ijkl}^3 (\mathbf{r})$ is a discontinuous scalar field (representing volumetric charge densities, for example).

In 2-D, as in the numerical examples considered here, $W_i^0 (\mathbf{r})$
and $\mathbf{W}_{ij}^1$ write exactly as above, but ${W}_{ijk}^2$ reduces to a scalar discontinuous function
\begin{flalign}
{W}_{ijk}^2 (\mathbf{r})=
	\begin{cases}
	\frac{1}{A}, \text{ if } \mathbf{r} \text{ is in the triangle $ijk$},
	\\ 0, \text{ otherwise},
	\end{cases}
\end{flalign}
where $A$ is the area of the triangle $ijk$~\footnote{Alternatively, one could consider it as a discontinuous vector function with such amplitude and oriented along the $z$-direction, i.e., transverse to a 2-D domain in the $xy$-plane, so that expressions such as~(\ref{hodge2}) remain invariant with the volume element $dV$ representing an area (2-D volume).}. Furthermore, ${W}_{ijkl}^3$ is identically zero in 2-D.
For these and more properties of Whitney forms, the reader is refered to~\cite{Bossavit88:Whitney, Bossavit02:Generating, Arnold06:Finite, Teixeira13:Differential, Teixeira14:Lattice} and references therein.

%%%%%%%%%%%%%%%%%%%%%%%%%%%%%%%%%%%%%%%%%%%%%%%%%%%%%%%%%%%%%%%%%%%%%%%%%%%%%%%%%%%%%%%%%%%%%%%%%%%%%%%%%%%%%%
\section{Line integral of Whitney 1-forms}
\label{app.b}
\begin{figure}[t]
    \centering
    \includegraphics[width=2.8in]{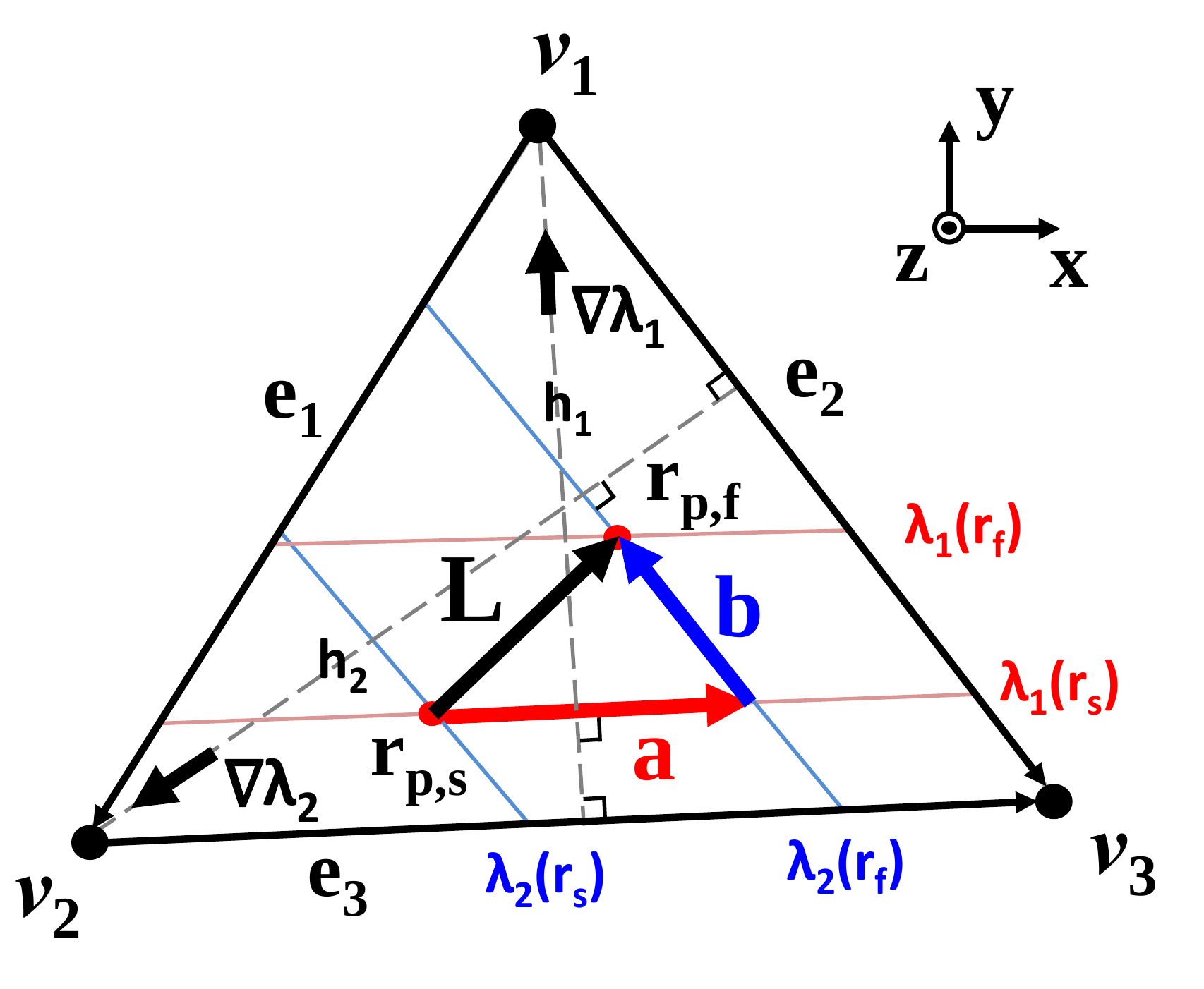}
    \caption{An arbitrary particle path $\mathbf{L}$ during $\Delta t$ and associated parameters in the $xy$-plane.}
    \label{arb.path}
\end{figure}

%%%%%%%%%%%%%%%%%%%%%%%%%%%%%%%%%%%%%%%%%%%%%%%%%%%%%%%%%%%%%%%%%%%%%%%%%%%%%%%%%%%%%%%%%%%%%%%%%%%%%%%%%%%%%%
The scatter step of the proposed algorithm and the analytical verification of charge conservation provided above both rely upon the evaluation of line integrals of Whitney 1-forms. In this Appendix, we consider this in more detail. An arbitrary segment $\mathbf{L}$ from $\mathbf{r}_{p,s}$ to $\mathbf{r}_{p,f}$ on a triangle is illustrated in Fig. \ref{arb.path}. The segment can be decomposed into two vectors $\mathbf{a}$ and $\mathbf{b}$. $\lambda_1(\cdot)$ and $\lambda_2(\cdot)$ are barycentric coordinates associated with $\nu_1$ and $\nu_2$. $h_1$ and $h_2$ are the heights of the triangle for the base of $e_3$ and $e_2$, respectively. The edge vectors $\mathbf{e}_1$, $\mathbf{e}_2$, and $\mathbf{e}_3$ are oriented in an ascending fashion of the associated vertex numbers. Note that the edge numbers do not coincide with the vertex numbers.

A simple way to evaluate the line integral
\begin{flalign} \label{line.Whitney}
\int_{\mathbf{r}_{p,s}}^{\mathbf{r}_{p,f}} \mathbf{W}_i^1(\mathbf{r}_p)\cdot d{\mathbf{L}},
\end{flalign}
is to use a parametric representation such that
\begin{flalign}
\mathbf{W}_i^1(\mathbf{r}_p)=\mathbf{W}_i^1(s)\quad \text{and}\quad d{\mathbf{L}}=d{\mathbf{L}}(s).
\end{flalign}
$\mathbf{W}_i^1(s)$ and $d{\mathbf{L}}(s)$ are simply assumed to be a linear function of the parameter $s$ and the range of $s$ is set to be $0\leq s \leq 1$. As an example, the line integral of the Whitney edge basis function associated with $e_1$, $\mathbf{W}_1^1(\mathbf{r}_p)=\lambda_1\na\lambda_2-\lambda_2\na\lambda_1$, is derived here.

As preliminaries, several variables are calculated. Vectors $\mathbf{a}$ and $\mathbf{b}$ can be expressed as
\begin{flalign}
\mathbf{a}=-\left(\lambda_2^f - \lambda_2^s\right)\mathbf{e}_3 = -\Delta_2 \mathbf{e}_3, \\
\mathbf{b}=-\left(\lambda_1^f - \lambda_1^s\right)\mathbf{e}_2 = -\Delta_1 \mathbf{e}_2.
\end{flalign}
First, the path and its space derivative are parameterized through $s$ such that
\begin{flalign}
\mathbf{L}&=\mathbf{L}(s)=(\mathbf{a}+\mathbf{b})s=-(\Delta_2 \mathbf{e}_3+\Delta_1 \mathbf{e}_2)s,
\end{flalign}
and $d\mathbf{L}=-(\Delta_2 \mathbf{e}_3+\Delta_1 \mathbf{e}_2)ds$.
Next, barycentric coordinates and their gradients are parameterized through $s$ as well, i.e.,
\begin{flalign}
\lambda_1(s)=\left[\lambda_1^f - \lambda_1^s\right]s+\lambda_1^s=\Delta_1 s + \lambda_1^s \\
\lambda_2(s)=\left[\lambda_2^f - \lambda_2^s\right]s+\lambda_2^s=\Delta_2 s + \lambda_2^s
\end{flalign}
The gradients of the barycentric coordinates are constant, so they are not the function of $s$, that is
\begin{flalign}
\na\lambda_1=\frac{1}{2A}\hat{z}\times\mathbf{e}_3, \quad
\na\lambda_2=\frac{1}{2A}\mathbf{e}_2\times\hat{z},
\end{flalign}
where $A$ is the area of the triangle. Some dot products used for the line integral are summarized below.
\begin{subequations}
\begin{flalign}
\na\lambda_1\cdot\mathbf{e}_2&=-1, \quad
\na\lambda_1\cdot\mathbf{e}_3=0 \\
\na\lambda_2\cdot\mathbf{e}_2&=0, \quad\;\;\;
\na\lambda_2\cdot\mathbf{e}_3=-1
\end{flalign}
\end{subequations}
Therefore, \eqref{line.Whitney} for $e_1$ is computed as
\allowdisplaybreaks
\begin{flalign}
&\int_{\mathbf{r}_{p,s}}^{\mathbf{r}_{p,f}} \mathbf{W}_1^1(\mathbf{r}_p)\cdot d{\mathbf{L}} =\int_{\mathbf{r}_{p,s}}^{\mathbf{r}_{p,f}}\left(\lambda_1\na\lambda_2-\lambda_2\na\lambda_1\right)
	\cdot d{\mathbf{L}} \notag\\
&=\int_0^1 \left[\left(\Delta_1 s + \lambda_1^s\right)\na\lambda_2
	- \left(\Delta_2 s + \lambda_2^s\right)\na\lambda_1\right]
	\cdot\left(-\Delta_2 \mathbf{e}_3 - \Delta_1 \mathbf{e}_2\right)ds \notag\\
&=-\Delta_2\left(\na\lambda_2\cdot\mathbf{e}_3\right) \int_0^1\left(\Delta_1 s + \lambda_1^s\right) ds
	+\Delta_1\left(\na\lambda_1\cdot\mathbf{e}_2\right) \int_0^1\left(\Delta_2 s + \lambda_2^s\right) ds \notag\\
&=\Delta_2\left[\frac{\Delta_1}{2}+\lambda_1^s\right]
	-\Delta_1\left[\frac{\Delta_2+\lambda_2^s}{2}\right] =\Delta_2\lambda_1^s - \Delta_1\lambda_2^s \notag\\
&=\left(\lambda_2^f-\lambda_2^s\right)\lambda_1^s - \left(\lambda_1^f-\lambda_1^s\right)\lambda_2^s =\lambda_1^s\lambda_2^f - \lambda_1^f\lambda_2^s.
\end{flalign}
Similarly, the other two line integrals can be computed as
\begin{flalign}
\int_{\mathbf{r}_{p,s}}^{\mathbf{r}_{p,f}} \mathbf{W}_2^1(\mathbf{r}_p)\cdot d{\mathbf{L}}
&=\lambda_1^s\lambda_3^f - \lambda_1^f\lambda_3^s
\end{flalign}
\begin{flalign}
\int_{\mathbf{r}_{p,s}}^{\mathbf{r}_{p,f}} \mathbf{W}_3^1(\mathbf{r}_p)\cdot d{\mathbf{L}}
&=\lambda_2^s\lambda_3^f - \lambda_2^f\lambda_3^s
\end{flalign}

%% References


\bibliographystyle{model1-num-names}
\bibliography{refsCPC}

\begin{thebibliography}{46}
\expandafter\ifx\csname natexlab\endcsname\relax\def\natexlab#1{#1}\fi
\providecommand{\bibinfo}[2]{#2}
\ifx\xfnm\relax \def\xfnm[#1]{\unskip,\space#1}\fi
%Type = Article
\bibitem[{Campos-Pinto et~al.(2014)Campos-Pinto, Jund, Salmon, and
  Sonnendr\"{u}cker}]{Pinto14:Charge}
\bibinfo{author}{M.~Campos-Pinto}, \bibinfo{author}{S.~Jund},
  \bibinfo{author}{S.~Salmon}, \bibinfo{author}{E.~Sonnendr\"{u}cker},
\newblock \bibinfo{title}{Charge-conserving {FEM}-{PIC} schemes on general
  grids},
\newblock \bibinfo{journal}{C. R. Mec.} \bibinfo{volume}{342}
  (\bibinfo{year}{2014}) \bibinfo{pages}{570--582}.
%Type = Article
\bibitem[{Squire et~al.(2012)Squire, Qin, and Tang}]{Squire12:Geometric}
\bibinfo{author}{J.~Squire}, \bibinfo{author}{H.~Qin}, \bibinfo{author}{W.~M.
  Tang},
\newblock \bibinfo{title}{Geometric integration of the {V}lasov-{M}axwell
  system with a variational particle-in-cell scheme},
\newblock \bibinfo{journal}{Phys. Plasmas} \bibinfo{volume}{19}
  (\bibinfo{year}{2012}).
%Type = Book
\bibitem[{Hockney and Eastwood(1981)}]{Hockney:Computer}
\bibinfo{author}{R.~W. Hockney}, \bibinfo{author}{J.~W. Eastwood},
  \bibinfo{title}{Computer Simulation Using Particles},
  \bibinfo{publisher}{McGraw-Hill}, \bibinfo{address}{New York},
  \bibinfo{year}{1981}.
%Type = Book
\bibitem[{Birdsall and Langdon(1985)}]{Birdsall:Plasma}
\bibinfo{author}{C.~K. Birdsall}, \bibinfo{author}{A.~B. Langdon},
  \bibinfo{title}{Plasma Physics via Computer Simulation},
  \bibinfo{publisher}{McGraw-Hill}, \bibinfo{address}{New York},
  \bibinfo{year}{1985}.
%Type = Book
\bibitem[{Fehske et~al.(2008)Fehske, Schneider, and
  Wei{\ss}e}]{Fehske:Computational}
\bibinfo{editor}{H.~Fehske}, \bibinfo{editor}{R.~Schneider},
  \bibinfo{editor}{A.~Wei{\ss}e} (Eds.), \bibinfo{title}{Computational
  Many-Particle Physics}, Lecture Notes in Physics,
  \bibinfo{publisher}{Springer}, \bibinfo{address}{Berlin, Germany},
  \bibinfo{year}{2008}.
%Type = Article
\bibitem[{Dawson(1983)}]{Dawson83:Particle}
\bibinfo{author}{J.~M. Dawson},
\newblock \bibinfo{title}{Particle simulation of plasmas},
\newblock \bibinfo{journal}{Rev. Mod. Phys.} \bibinfo{volume}{55}
  (\bibinfo{year}{1983}) \bibinfo{pages}{403--447}.
%Type = Article
\bibitem[{Bruhwiler et~al.(2001)Bruhwiler, Giacone, Cary, Verboncoeur, Mardahl,
  Esarey, Leemans, and Shadwick}]{Bruhwiler01:Particle}
\bibinfo{author}{D.~L. Bruhwiler}, \bibinfo{author}{R.~E. Giacone},
  \bibinfo{author}{J.~R. Cary}, \bibinfo{author}{J.~P. Verboncoeur},
  \bibinfo{author}{P.~Mardahl}, \bibinfo{author}{E.~Esarey},
  \bibinfo{author}{W.~P. Leemans}, \bibinfo{author}{B.~A. Shadwick},
\newblock \bibinfo{title}{Particle-in-cell simulations of plasma accelerators
  and electron-neutral collisions},
\newblock \bibinfo{journal}{Phys. Rev. Spec. Top. Accel. Beams}
  \bibinfo{volume}{4} (\bibinfo{year}{2001}).
%Type = Article
\bibitem[{Strozzi et~al.(2012)Strozzi, Tabak, Larson, Divol, Kemp, Bellei,
  Marinak, and Key}]{Strozzi12:Fast}
\bibinfo{author}{D.~J. Strozzi}, \bibinfo{author}{M.~Tabak},
  \bibinfo{author}{D.~J. Larson}, \bibinfo{author}{L.~Divol},
  \bibinfo{author}{A.~J. Kemp}, \bibinfo{author}{C.~Bellei},
  \bibinfo{author}{M.~M. Marinak}, \bibinfo{author}{M.~H. Key},
\newblock \bibinfo{title}{Fast-ignition transport studies: {R}ealistic electron
  source, integrated particle-in-cell and hydrodynamic modeling, imposed
  magnetic fields},
\newblock \bibinfo{journal}{Phys. Plasmas} \bibinfo{volume}{19}
  (\bibinfo{year}{2012}).
%Type = Article
\bibitem[{Booske(2008)}]{Booske08:Plasma}
\bibinfo{author}{J.~H. Booske},
\newblock \bibinfo{title}{Plasma physics and related challenges of
  millimeter-wave-to-terahertz and high power microwave generation},
\newblock \bibinfo{journal}{Phys. Plasmas} \bibinfo{volume}{15}
  (\bibinfo{year}{2008}).
%Type = Article
\bibitem[{Marder(1987)}]{Marder87:Method}
\bibinfo{author}{B.~Marder},
\newblock \bibinfo{title}{A method for incorporating {G}auss' law into
  electromagnetic {PIC} codes},
\newblock \bibinfo{journal}{J. Comput. Phys.} \bibinfo{volume}{68}
  (\bibinfo{year}{1987}) \bibinfo{pages}{48--55}.
%Type = Article
\bibitem[{Langdon(1992)}]{Langdon92:Enforcing}
\bibinfo{author}{A.~B. Langdon},
\newblock \bibinfo{title}{On enforcing {G}auss' law in electromagnetic
  particle-in-cell codes},
\newblock \bibinfo{journal}{Comput. Phys. Commun.} \bibinfo{volume}{70}
  (\bibinfo{year}{1992}) \bibinfo{pages}{447--450}.
%Type = Article
\bibitem[{Mardahl and Verboncoeur(1997)}]{Mardahl97:Charge}
\bibinfo{author}{P.~J. Mardahl}, \bibinfo{author}{J.~P. Verboncoeur},
\newblock \bibinfo{title}{Charge conservation in electromagnetic {PIC} codes;
  spectral comparison of {B}oris/{DADI} and {L}angdon-{M}arder methods},
\newblock \bibinfo{journal}{Comput. Phys. Commun.} \bibinfo{volume}{106}
  (\bibinfo{year}{1997}) \bibinfo{pages}{219--229}.
%Type = Article
\bibitem[{Eastwood(1991)}]{Eastwood91:Virtual}
\bibinfo{author}{J.~W. Eastwood},
\newblock \bibinfo{title}{The virtual particle electromagnetic particle-mesh
  method},
\newblock \bibinfo{journal}{Comput. Phys. Commun.} \bibinfo{volume}{64}
  (\bibinfo{year}{1991}) \bibinfo{pages}{252--266}.
%Type = Article
\bibitem[{Villasenor and Buneman(1992)}]{Villasenor92:Rigorous}
\bibinfo{author}{J.~Villasenor}, \bibinfo{author}{O.~Buneman},
\newblock \bibinfo{title}{Rigorous charge conservation for local
  electromagnetic field solvers},
\newblock \bibinfo{journal}{Comput. Phys. Commun.} \bibinfo{volume}{69}
  (\bibinfo{year}{1992}) \bibinfo{pages}{306--316}.
%Type = Article
\bibitem[{Esirkepov(2001)}]{Esirkepov01:Exact}
\bibinfo{author}{T.~Z. Esirkepov},
\newblock \bibinfo{title}{Exact charge conservation scheme for
  {P}article-in-{C}ell simulation with an arbitrary form-factor},
\newblock \bibinfo{journal}{Comput. Phys. Commun.} \bibinfo{volume}{135}
  (\bibinfo{year}{2001}) \bibinfo{pages}{144--153}.
%Type = Article
\bibitem[{Umeda et~al.(2003)Umeda, Omura, Tominaga, and
  Matsumoto}]{Umeda03:New}
\bibinfo{author}{T.~Umeda}, \bibinfo{author}{Y.~Omura},
  \bibinfo{author}{T.~Tominaga}, \bibinfo{author}{H.~Matsumoto},
\newblock \bibinfo{title}{A new charge conservation method in electromagnetic
  particle-in-cell simulations},
\newblock \bibinfo{journal}{Comput. Phys. Commun.} \bibinfo{volume}{156}
  (\bibinfo{year}{2003}) \bibinfo{pages}{73--85}.
%Type = Article
\bibitem[{Kong et~al.(2011)Kong, Huang, Ren, and Decyk}]{Kong11:Particle}
\bibinfo{author}{X.~Kong}, \bibinfo{author}{M.~C. Huang},
  \bibinfo{author}{C.~Ren}, \bibinfo{author}{V.~K. Decyk},
\newblock \bibinfo{title}{Particle-in-cell simulations with charge-conserving
  current deposition on graphic processing units},
\newblock \bibinfo{journal}{J. Comput. Phys.} \bibinfo{volume}{230}
  (\bibinfo{year}{2011}) \bibinfo{pages}{1676--1685}.
%Type = Article
\bibitem[{Sokolov(2013)}]{Sokolov13:Alternating}
\bibinfo{author}{I.~V. Sokolov},
\newblock \bibinfo{title}{Alternating-order interpolation in a
  charge-conserving scheme for particle-in-cell simulations},
\newblock \bibinfo{journal}{Comput. Phys. Commun.} \bibinfo{volume}{184}
  (\bibinfo{year}{2013}) \bibinfo{pages}{320--328}.
%Type = Article
\bibitem[{Jacobs and Hesthaven(2006)}]{Jacobs06:High}
\bibinfo{author}{G.~B. Jacobs}, \bibinfo{author}{J.~S. Hesthaven},
\newblock \bibinfo{title}{High-order nodal discontinuous {G}alerkin
  particle-in-cell method on unstructured grids},
\newblock \bibinfo{journal}{J. Comput. Phys.} \bibinfo{volume}{214}
  (\bibinfo{year}{2006}) \bibinfo{pages}{96--121}.
%Type = Inproceedings
\bibitem[{Candel et~al.(2007)Candel, Kabel, Lee, Li, Limborg, Ng, Prudencio,
  Schussman, Uplenchwar, and Ko}]{Candel07:Parallel}
\bibinfo{author}{A.~Candel}, \bibinfo{author}{A.~Kabel},
  \bibinfo{author}{L.~Lee}, \bibinfo{author}{Z.~Li},
  \bibinfo{author}{C.~Limborg}, \bibinfo{author}{C.~Ng},
  \bibinfo{author}{E.~Prudencio}, \bibinfo{author}{G.~Schussman},
  \bibinfo{author}{R.~Uplenchwar}, \bibinfo{author}{K.~Ko},
\newblock \bibinfo{title}{Parallel finite element particle-in-cell code for
  simulations of space-charge dominated beam-cavity interactions},
\newblock in: \bibinfo{booktitle}{Proc. IEEE Particle Accelerator
  Conference'07}, \bibinfo{address}{Albuquerque, NM}, pp.
  \bibinfo{pages}{908--910}.
%Type = Techreport
\bibitem[{Candel et~al.(2009)Candel, Kabel, Lee, Li, Limborg, Ng, Prudencio,
  Schussman, Uplenchwar, and Ko}]{Candel09:Parallel}
\bibinfo{author}{A.~Candel}, \bibinfo{author}{A.~Kabel},
  \bibinfo{author}{L.~Lee}, \bibinfo{author}{Z.~Li},
  \bibinfo{author}{C.~Limborg}, \bibinfo{author}{C.~Ng},
  \bibinfo{author}{E.~Prudencio}, \bibinfo{author}{G.~Schussman},
  \bibinfo{author}{R.~Uplenchwar}, \bibinfo{author}{K.~Ko},
  \bibinfo{title}{Parallel higher-order finite element method for accurate
  field computations in wakefield and {PIC} simulations}, \bibinfo{type}{Tech.
  Rep.} \bibinfo{number}{{SLAC}-{PUB}-13667}, SLAC, \bibinfo{address}{Menlo
  Park, CA}, \bibinfo{year}{2009}.
%Type = Article
\bibitem[{Manges and Cendes(1995)}]{Manges95:cotree}
\bibinfo{author}{J.~B. Manges}, \bibinfo{author}{Z.~J. Cendes},
\newblock \bibinfo{title}{A generalized tree-cotree gauge for magnetic field
  computation},
\newblock \bibinfo{journal}{IEEE Trans. Magn.} \bibinfo{volume}{31}
  (\bibinfo{year}{1995}) \bibinfo{pages}{1342--1347}.
%Type = Article
\bibitem[{Albanese and Rubinacci(1988)}]{Albanese88:eddy}
\bibinfo{author}{R.~Albanese}, \bibinfo{author}{G.~Rubinacci},
\newblock \bibinfo{title}{Integral formulation for 3{D} eddy-current
  computation using edge elements},
\newblock \bibinfo{journal}{Proc. IEE pt.A} \bibinfo{volume}{135}
  (\bibinfo{year}{1988}) \bibinfo{pages}{457--462}.
%Type = Article
\bibitem[{Hwang and Wu(1999)}]{Hwang99:stability}
\bibinfo{author}{C.-T. Hwang}, \bibinfo{author}{R.-B. Wu},
\newblock \bibinfo{title}{Treating late-time instability of hybrid
  finite-element/finite-difference time-domain method},
\newblock \bibinfo{journal}{IEEE Trans. Antennas Propag.} \bibinfo{volume}{47}
  (\bibinfo{year}{1999}) \bibinfo{pages}{227--232}.
%Type = Article
\bibitem[{Moon et~al.(2014)Moon, Teixeira, Kim, and Omelchenko}]{Moon14:Trade}
\bibinfo{author}{H.~Moon}, \bibinfo{author}{F.~L. Teixeira},
  \bibinfo{author}{J.~Kim}, \bibinfo{author}{Y.~A. Omelchenko},
\newblock \bibinfo{title}{Trade-offs for unconditional stability in the
  finite-element time-domain method},
\newblock \bibinfo{journal}{IEEE Microw. Wireless Compon. Lett.}
  \bibinfo{volume}{24} (\bibinfo{year}{2014}) \bibinfo{pages}{361--363}.
%Type = Article
\bibitem[{Kim and Teixeira(2011)}]{Kim11:Parallel}
\bibinfo{author}{J.~Kim}, \bibinfo{author}{F.~L. Teixeira},
\newblock \bibinfo{title}{Parallel and explicit finite-element time-domain
  method for {M}axwell's equations},
\newblock \bibinfo{journal}{IEEE Trans. Antennas Propag.} \bibinfo{volume}{59}
  (\bibinfo{year}{2011}) \bibinfo{pages}{2350--2356}.
%Type = Article
\bibitem[{He and Teixeira(2006)}]{He06:Geometric}
\bibinfo{author}{B.~He}, \bibinfo{author}{F.~L. Teixeira},
\newblock \bibinfo{title}{Geometric finite element discretization of {M}axwell
  equations in primal and dual spaces},
\newblock \bibinfo{journal}{Phys. Lett. A} \bibinfo{volume}{349}
  (\bibinfo{year}{2006}) \bibinfo{pages}{1--14}.
%Type = Article
\bibitem[{Bossavit(1988)}]{Bossavit88:Whitney}
\bibinfo{author}{A.~Bossavit},
\newblock \bibinfo{title}{Whitney forms: a class of finite-elements for
  three-dimensional computations in electromagnetism},
\newblock \bibinfo{journal}{IEE Proc. A} \bibinfo{volume}{135}
  (\bibinfo{year}{1988}) \bibinfo{pages}{493--500}.
%Type = Article
\bibitem[{Bossavit(2002)}]{Bossavit02:Generating}
\bibinfo{author}{A.~Bossavit},
\newblock \bibinfo{title}{Generating {W}hitney forms of polynomial degree one
  and higher},
\newblock \bibinfo{journal}{IEEE Trans. Magn.} \bibinfo{volume}{38}
  (\bibinfo{year}{2002}) \bibinfo{pages}{341--344}.
%Type = Article
\bibitem[{Deschamps(1981)}]{Deschamps81:Electromagnetics}
\bibinfo{author}{G.~A. Deschamps},
\newblock \bibinfo{title}{Electromagnetics and differential forms},
\newblock \bibinfo{journal}{Proc. of the IEEE} \bibinfo{volume}{69}
  (\bibinfo{year}{1981}) \bibinfo{pages}{676--696}.
%Type = Article
\bibitem[{Warnick et~al.(1997)Warnick, Selfridge, and
  Arnold}]{Warnick97:Teaching}
\bibinfo{author}{K.~F. Warnick}, \bibinfo{author}{R.~H. Selfridge},
  \bibinfo{author}{D.~V. Arnold},
\newblock \bibinfo{title}{Teaching electromagnetic field theory using
  differential forms},
\newblock \bibinfo{journal}{IEEE Trans. Educ.} \bibinfo{volume}{40}
  (\bibinfo{year}{1997}) \bibinfo{pages}{53--68}.
%Type = Article
\bibitem[{Teixeira and Chew(1999)}]{Teixeira99:Lattice}
\bibinfo{author}{F.~L. Teixeira}, \bibinfo{author}{W.~C. Chew},
\newblock \bibinfo{title}{Lattice electromagnetic theory from a topological
  viewpoint},
\newblock \bibinfo{journal}{J. Math. Phys.} \bibinfo{volume}{40}
  (\bibinfo{year}{1999}) \bibinfo{pages}{169--187}.
%Type = Article
\bibitem[{He and Teixeira(2007)}]{He07:Differential}
\bibinfo{author}{B.~He}, \bibinfo{author}{F.~L. Teixeira},
\newblock \bibinfo{title}{Differential forms, {G}alerkin duality, and sparse
  inverse approximations in finite element solutions of {M}axwell equations},
\newblock \bibinfo{journal}{IEEE Trans. Antennas Propag.} \bibinfo{volume}{55}
  (\bibinfo{year}{2007}) \bibinfo{pages}{1359--1368}.
%Type = Article
\bibitem[{Clemens and Weiland(2001)}]{Clemens01:Discrete}
\bibinfo{author}{M.~Clemens}, \bibinfo{author}{T.~Weiland},
\newblock \bibinfo{title}{Discrete electromagnetism with the finite integration
  technique},
\newblock \bibinfo{journal}{Prog. Electromagn. Res.} \bibinfo{volume}{32}
  (\bibinfo{year}{2001}) \bibinfo{pages}{65--87}.
%Type = Article
\bibitem[{Schuhmann and Weiland(2001)}]{Schuhmann01:Conservation}
\bibinfo{author}{R.~Schuhmann}, \bibinfo{author}{T.~Weiland},
\newblock \bibinfo{title}{Conservation of discrete energy and related laws in
  the finite integration technique},
\newblock \bibinfo{journal}{Prog. Electromagn. Res.} \bibinfo{volume}{32}
  (\bibinfo{year}{2001}) \bibinfo{pages}{301--316}.
%Type = Article
\bibitem[{Lee(2006)}]{Lee06:Note}
\bibinfo{author}{R.~Lee},
\newblock \bibinfo{title}{A note on mass lumping in the finite element time
  domain method},
\newblock \bibinfo{journal}{IEEE Trans. Antennas Propag.} \bibinfo{volume}{54}
  (\bibinfo{year}{2006}) \bibinfo{pages}{760--762}.
%Type = Article
\bibitem[{He and Teixeira(2005)}]{He05:Degrees}
\bibinfo{author}{B.~He}, \bibinfo{author}{F.~L. Teixeira},
\newblock \bibinfo{title}{On the degrees of freedom of lattice
  electrodynamics},
\newblock \bibinfo{journal}{Phys. Lett. A} \bibinfo{volume}{336}
  (\bibinfo{year}{2005}) \bibinfo{pages}{1--7}.
%Type = Book
\bibitem[{Flanders(1989)}]{Flanders:Differential}
\bibinfo{author}{H.~Flanders}, \bibinfo{title}{Differential Forms with
  Applications to the Physical Sciences}, \bibinfo{publisher}{Dover
  Publications}, \bibinfo{address}{Mineola, N.Y.}, \bibinfo{year}{1989}.
%Type = Book
\bibitem[{Whitney(1957)}]{Whitney:Geometric}
\bibinfo{author}{H.~Whitney}, \bibinfo{title}{Geometric Integration Theory},
  Princeton Mathematical Series, \bibinfo{publisher}{Princeton University
  Press}, \bibinfo{address}{Princeton}, \bibinfo{year}{1957}.
%Type = Book
\bibitem[{Bossavit(1998)}]{Bossavit:Computational}
\bibinfo{author}{A.~Bossavit}, \bibinfo{title}{Computational Electromagnetism:
  Variational Formulations, Complementarity, Edge Elements}, Electromagnetism,
  \bibinfo{publisher}{Academic Press}, \bibinfo{address}{San Diego},
  \bibinfo{year}{1998}.
%Type = Book
\bibitem[{Jin(2002)}]{Jin:Finite}
\bibinfo{author}{J.-M. Jin}, \bibinfo{title}{The Finite Element Method in
  Electromagnetics}, \bibinfo{publisher}{Wiley}, \bibinfo{address}{New York},
  \bibinfo{edition}{second} edition, \bibinfo{year}{2002}.
%Type = Book
\bibitem[{Bondeson et~al.(2005)Bondeson, Rylander, and
  Ingelstr{\"o}m}]{Bondeson:Computational}
\bibinfo{author}{A.~Bondeson}, \bibinfo{author}{T.~Rylander},
  \bibinfo{author}{P.~Ingelstr{\"o}m}, \bibinfo{title}{Computational
  Electromagnetics}, Texts in Applied Mathematics,
  \bibinfo{publisher}{Springer}, \bibinfo{address}{New York, N.Y.},
  \bibinfo{year}{2005}.
%Type = Book
\bibitem[{Silvester and Ferrari(1996)}]{Silvester:Finite}
\bibinfo{author}{P.~P. Silvester}, \bibinfo{author}{R.~L. Ferrari},
  \bibinfo{title}{Finite Elements for Electrical Engineers},
  \bibinfo{publisher}{Cambridge University Press}, \bibinfo{address}{New York},
  \bibinfo{edition}{third} edition, \bibinfo{year}{1996}.
%Type = Article
\bibitem[{Teixeira(2014)}]{Teixeira14:Lattice}
\bibinfo{author}{F.~L. Teixeira},
\newblock \bibinfo{title}{Lattice {M}axwell's equations},
\newblock \bibinfo{journal}{Prog. Electromagn. Res.} \bibinfo{volume}{148}
  (\bibinfo{year}{2014}) \bibinfo{pages}{113--128}.
%Type = Article
\bibitem[{Arnold et~al.(2006)Arnold, Falk, and Winther}]{Arnold06:Finite}
\bibinfo{author}{D.~N. Arnold}, \bibinfo{author}{R.~S. Falk},
  \bibinfo{author}{R.~Winther},
\newblock \bibinfo{title}{Finite element exterior calculus, homological
  techniques, and applications},
\newblock \bibinfo{journal}{Acta Numerica} \bibinfo{volume}{15}
  (\bibinfo{year}{2006}) \bibinfo{pages}{1--155}.
%Type = Article
\bibitem[{Teixeira(2013)}]{Teixeira13:Differential}
\bibinfo{author}{F.~L. Teixeira},
\newblock \bibinfo{title}{Differential forms in lattice field theories: {A}n
  overview},
\newblock \bibinfo{journal}{ISRN Math. Phys.} \bibinfo{volume}{2013}
  (\bibinfo{year}{2013}) \bibinfo{pages}{16}.

\end{thebibliography}






% \end{thebibliography}

\end{document}